\documentclass[11pt]{article}

\usepackage[a4paper,margin=1in]{geometry}
\usepackage[T1]{fontenc}
\usepackage{lmodern}
\usepackage{microtype}
\usepackage{amsmath,amssymb,amsthm}
\usepackage{aliascnt}
\usepackage{tikz-cd}
\usepackage{enumitem}
\usepackage[hidelinks]{hyperref}
\usepackage[nameinlink,noabbrev]{cleveref}

\setlist[itemize]{leftmargin=1.5em}
\setlist[enumerate]{leftmargin=1.7em}
\newtheorem{theorem}{Theorem}[section]
\newaliascnt{proposition}{theorem}
\newtheorem{proposition}[proposition]{Proposition}
\aliascntresetthe{proposition}
\newaliascnt{lemma}{theorem}
\newtheorem{lemma}[lemma]{Lemma}
\aliascntresetthe{lemma}
\newaliascnt{corollary}{theorem}
\newtheorem{corollary}[corollary]{Corollary}
\aliascntresetthe{corollary}

\theoremstyle{definition}
\newaliascnt{definition}{theorem}
\newtheorem{definition}[definition]{Definition}
\aliascntresetthe{definition}

\theoremstyle{remark}
\newaliascnt{remark}{theorem}
\newtheorem{remark}[remark]{Remark}
\aliascntresetthe{remark}
\newtheorem*{claim}{Claim}
\newtheorem*{researchquestion}{Question}

\newcommand{\F}{\mathbb{F}}

\newcommand{\cl}{\operatorname{cl}}
\newcommand{\crat}[5]{\left[\begin{matrix}#1&#2\\#3&#4\end{matrix}\right]_{#5}}

\title{Five-Entropic Matroids Are Quinary}

\author{
	\begin{tabular}{c@{\hspace{.5cm}}c}
		Mohammad Hossein Kalantari\thanks{\href{mailto:mohammadhossein.kalantari12@gmail.com}{mohammadhossein.kalantari12@gmail.com}}
		&
		Shahram Khazaei\thanks{\href{mailto:shahram.khazaei@sharif.ir}{shahram.khazaei@sharif.ir}}
		\\[2mm]
		\multicolumn{2}{c}{\small Department of Mathematical Sciences, Sharif University of Technology}
		\\
		\multicolumn{2}{c}{\small Tehran, Iran}
	\end{tabular}
}

\date{}

\begin{document}
\maketitle
\begin{abstract}
A finite matroid is representable over a field if its elements
can be assigned vectors over that field such that the rank of
every subset equals the dimension of the span of its vectors.
For an integer $q\geq 2$, a matroid is $q$-entropic if there
exists a random vector indexed by its ground set, with each
coordinate taking values in a common $q$-element alphabet,
such that the entropy of every marginal subvector equals the
rank of the corresponding subset multiplied by $\log q$.

It remains open to determine for which values of $q$ every
$q$-entropic matroid is representable over a field.
For $q=2$ and $q=3$, representability over $\F_q$ has
long been known. More recently, it was established that
every $4$-entropic matroid is representable over $\F_4$
and that the $6$-entropic matroids are precisely the
regular matroids, that is, those representable over
every field.

We settle the previously unresolved five-symbol case,
proving that a finite matroid is $5$-entropic if and
only if it is representable over $\F_5$.
For the nontrivial implication, we use Mat\'u\v{s}'s
equivalence between entropic and partition representations.
We show that a degree-five partition representation
induces, on each embedded $U_{2,4}$-minor, an
$\F_5$-linear representation after coordinatewise
relabeling, with a rescaling class independent of the
chosen linearization. Using the foundation presentation
of Baker and Lorscheid, we verify the compatibility
relations among these local classes and thereby obtain
a global $\F_5$-representation.
\end{abstract}


\section{Introduction}

Matroids abstract the notion of dependence by assigning
a rank to each subset of a ground set.
Their rank functions admit several notions of representation.
Entropic matroids admit equivalent descriptions in terms of
almost-affine codes~\cite{SimonisAshikhmin1998},
partition representations~\cite{Matus1999},
and ideal secret-sharing schemes
~\cite{BrickellDavenport1991,AbbeSpirkl2019}.
This paper proves that a finite matroid is $5$-entropic
if and only if it is representable over $\F_5$.

\subsection{Matroid representability and problem statement}

Let $M=(E,r)$ be a finite matroid. A representation of
$M$ over a field $K$ is a family of vectors
$(v_e)_{e\in E}$ such that
\[
r(A)=\dim_K\operatorname{span}\{v_e:e\in A\}
\qquad (A\subseteq E).
\]
We say that $M$ is $K$-representable if it admits
such a representation, and simply representable if
it is representable over some field.

A more general notion replaces vectors by subspaces.
For a positive integer $t$, a $t$-folded representation
of $M$ over $K$ is a family of subspaces $(V_e)_{e\in E}$
of a common vector space over $K$ such that
\[
\dim_K\left(\sum_{e\in A}V_e\right)=t\,r(A)
\qquad (A\subseteq E).
\]
A matroid is folded-representable if it admits such a
representation for some field $K$ and positive
integer $t$. Ordinary representability is
the special case $t=1$.

For an integer $q\geq 2$, the matroid $M$ is
$q$-entropic if there exists a random vector
$(X_e)_{e\in E}$, with each coordinate taking values
in a common $q$-element alphabet, such that
\[
H(X_A)=r(A)\log q
\qquad (A\subseteq E).
\]
A matroid is called entropic if it is $q$-entropic
for some integer $q\geq 2$.

Every $t$-folded representation over a finite
field of order $s$ yields an $s^t$-entropic
representation. In particular, for every prime
power $q$, every $\F_q$-representable matroid
is $q$-entropic.

\paragraph{Problem statement.}
Whether every entropic matroid is folded-representable
remains open. Simonis and Ashikhmin raised this
question in their study of almost-affine codes and
ideal secret sharing~\cite{SimonisAshikhmin1998}.
A related, more specific fixed-alphabet problem asks:

\begin{researchquestion}
For which integers $q\geq 2$ is every $q$-entropic
matroid representable over some field?
\end{researchquestion}

When $q$ is a prime power, one can ask the stronger
question of representability over $\F_q$.

\subsection{Related work and the small-alphabet landscape}

The cases $q=2$ and $q=3$ were settled by Beimel
and Chor in their study of ideal secret-sharing schemes
~\cite{BeimelChor1994}. Their results imply that
every $q$-entropic matroid is representable over
$\F_q$ for $q\in\{2,3\}$.
Alternative matroid-theoretic proofs were later given
by Abbe and Spirkl~\cite{AbbeSpirkl2019}.

It is classical that the non-Pappus matroid is not
representable over any field, by Pappus's theorem
(see e.g., ~\cite{Oxley}).
Simonis and Ashikhmin nevertheless exhibited a
$2$-folded representation of this matroid over
$\F_3$, showing that it is $9$-entropic
~\cite{SimonisAshikhmin1998}.

The four-symbol case was recently resolved independently
by two different approaches. Kalantari and Khazaei
proved that every $4$-entropic matroid is representable
over $\F_4$, using minor closure, the excluded-minor
characterization of quaternary matroids, and properties
of quasigroups of order four
~\cite{KalantariKhazaei2026}.
Independently, Beimel, Ben-Efraim, Farr\`as, and Moya
established the same result, formulated in the language
of ideal secret sharing, using matroid-theoretic and
Latin-square methods
~\cite{BeimelBenEfraimFarrasMoya2026}.
The two proofs share some combinatorial ingredients
but differ in their overall organization.

Beimel and coauthors also characterized ideal
secret-sharing schemes over six-element domains.
In matroid language, their result implies that
a matroid is $6$-entropic if and only if it is
regular, that is, representable over every field
~\cite{BeimelBenEfraimFarrasMoya2026}.

\subsection{Main result and proof strategy}

The preceding results leave open whether every
$5$-entropic matroid is representable over $\F_5$.
Our main result answers this question affirmatively
and gives an exact characterization.

\begin{theorem}[Five-symbol rigidity]\label{thm:main}
For every finite matroid $M$,
\[
M\text{ is $5$-entropic}
\quad\Longleftrightarrow\quad
M\text{ is representable over }\F_5.
\]
\end{theorem}

The implication from $\F_5$-representability to
$5$-entropicity follows from the standard linear
entropy construction~\cite[Lemma~3]{AbbeSpirkl2019}.
The main task is therefore to prove the converse.

In order to show the other implication, consider a matroid $M = (E,r)$ that is $5$-entropic. We prove that $M$ is
representable over $\F_5$ by taking the following steps:
\begin{enumerate}
\item
By \Cref{prop:entropy-partition}, $M$ has a partition representation of degree five on the ground set 
$\Omega$.
\item
By \Cref{remark: coordinate partition}, we can choose the ground set $\Omega$ so that $\Omega \subseteq \F_5^E$ and $(0)_{e \in E} \in \Omega$.
\item
If $N$ is an embedded minor of $M$ isomorphic to $U_{2,4}$, then we can use \Cref{lem:minor-closure} to 
get a partition representation of degree five for $N$ on the ground set $\Omega_N$ from $\Omega$.
\item
By \Cref{lem: U24 linearity}, there are permutations $f_i$ on the set $\F_5$, for $1 \leq i \leq 4$, so 
that transforming $\Omega_N$ by $f = (f_1,f_2,f_3,f_4)$ gives a linear partition representation, i.e.,
there exists a matrix $A \in \F_5^{2 \times N}$ such that 
\[
f(\Omega_N) = \{xA \in \F_5^N: x \in \F_5^2\}.
\] 
By \Cref{prop: linearity}, $A$ is a representation for $N$ over $\F_5$. By \Cref{prop: A GP function iff representable over F}, it is possible to construct a Grassmann--Pl\"ucker function for $N$ with values in $\F_5$ from $A$. Also, see
\Cref{rem: equivalence of represenation}. (Here, we use $N$ to denote both the embedded minor and also
the ground set of the embedded minor. Occasionally, we will use the matroid to denote its ground set, but the intended meaning will be clear from the context.)
\item
By \Cref{lem: U24 linearity equivalence}, different Grassmann--Pl\"ucker functions obtained in the previous
step are all rescaling equivalent. By \Cref{thm: a morphism if and only if a GP function}, these rescaling
equivalent Grassmann--Pl\"ucker functions induce a pasture morphism $\Phi_N$ from the foundation $F_N$ of $N$ to $\F_5$.
\item
Up to this point, we have been able to associate a morphism $\Phi_N: F_N \rightarrow \F_5$ to every embedded minor $N$ of $M$ isomorphic to $U_{2,4}$. Using the universal property in \Cref{prop: universal property of tensor product}, it is possible to put these morphisms together to get a morphism 
\[
\Phi : \left(\bigotimes_{N}F_N\right) \rightarrow \F_5
\]
Here, $N$ ranges over all embedded minors $N$ of $M$ isomorphic to $U_{2,4}$.
\item
$M$ has a partition representation of degree five. Using \Cref{lem:minor-closure}, \Cref{thm: fano}, and 
\Cref{thm:dual fano}, it follows that $M$ does not have an embedded minor isomorphic to the Fano matroid 
$F_7$ or its dual $F_7^\ast$.
\item
In \Cref{sec:proofmain}, we show that the morphism $\Phi$ satisfies relations (R1), (R2), and (R3) in 
\Cref{thm: BL foundation representation}. Moreover, as shown in the previous step, there is no need to check (R4) relation. Therefore, by \Cref{prop: morphism factors through quotient}, there is a morphism
\[
\tilde{\Phi} : \left(\bigotimes_{N}F_N\right)//S \rightarrow \F_5.
\]
By \Cref{thm: BL foundation representation}, it follows that there is a morphism from the foundation $F_M$
of $M$ to $\F_5$. By \Cref{thm: a morphism if and only if a GP function} and \Cref{prop: A GP function iff representable over F}, $M$ is representable over $\F_5$.
\end{enumerate}

\section{Preliminaries}\label{sec:partition}

We assume familiarity with basic matroid theory and refer
the reader to Oxley~\cite{Oxley} for standard definitions.
We recall the notation needed below and introduce partition
representations in the sense of Mat\'u\v{s}~\cite{Matus1999},
the quasigroup structures arising from matroid circuits,
and the foundation of a matroid.

\subsection{Matroid notation and embedded minors}

Let $M$ be a matroid on a finite ground set $E$.
We write $r_M$, or simply $r$, for its rank function.
For $I,J\subseteq E$, deletion and contraction are denoted
by $M\backslash I$ and $M/J$, respectively.
The uniform matroid $U_{r,n}$ has rank $r$ on an
$n$-element ground set, with independent sets precisely
the subsets of size at most $r$.

An \emph{embedded minor} of $M$ is a minor of the form
$M\backslash I/J$, where $I$ and $J$ are disjoint,
$I$ is coindependent in $M$, and $J$ is independent
in $M$. We retain the specified sets $I$ and $J$
when referring to an embedded minor.

\subsection{Partition representations and quasigroups}

A partition representation encodes matroid ranks through
common refinements of partitions of a finite set.
We first recall this representation and then explain
how matroid circuits give rise to quasigroup operations
on the corresponding coordinate alphabets.

\begin{definition}[$q$-entropic matroids]\label{def:qentropic}
Let $q\ge 2$. A matroid $M=(E,r)$ is \emph{$q$-entropic} if there are random variables $(X_e)_{e\in E}$, each taking values in a common $q$-element alphabet, such that
\[
H(X_A)=r(A)\log q\qquad\text{for every }A\subseteq E.
\]
\end{definition}

\begin{definition}[Partition representations~\cite{Matus1999}]\label{def:partition}
A \emph{partition (p-)representation of degree $q$} of $M=(E,r)$ consists of a finite set $\Omega$ with $|\Omega|=q^{r(E)}$ and partitions $\xi_e$ of $\Omega$, $e\in E$, i.e., $\xi_e$ is a partition of the ground set $\Omega$ for each $e \in E$, such that for every $A\subseteq E$ the common refinement $\bigwedge_{e\in A}\xi_e$ has exactly $q^{r(A)}$ blocks, all of cardinality $q^{r(E)-r(A)}$.
\end{definition}

\begin{definition}[Quasigroup sets~\cite{Matus1999}]\label{def: quasigroup set}
Consider  non-empty sets $A,Q$. Moreover, suppose that $2\leq|A|$. A set $E \subseteq Q^A$ is called a quasigroup set if for every $i \in A$ and $x_{A-i} \in Q^{A-i}$, there is a unique $x_i \in Q$ such that $x_{A} \in E$. $x_A$ is the element of $Q^A$ composed of $x_i$ and $x_{A-i}$, i.e., its i-th coordinate is $x_i$ and its other coordinates are $x_{A-i}$.
\end{definition}

\begin{definition}[Quasigroup maps]\label{def: quasigroup map}
Consider non-empty sets $A,Q$. A map $f:Q^A \rightarrow Q$ is called a quasigroup if for every $i \in A$ and $x_{A-i} \in Q^{A-i}$, the map $f_{i,x_{A-i}}: Q \rightarrow Q$ that is obtained by fixing all variables except $i$ by $x_{A-i}$, is a bijection.
\end{definition}

\begin{remark}\label{remark: quasi map eq set}
Consider a quasigroup $f:Q^A \rightarrow Q$. The set $\{(x,f(x))\in Q^{A} \times Q:x \in Q^A\}$ is a quasigroup set. Conversely, consider a quasigroup set $E \subseteq Q^A$ and $i \in A$. By definition of quasigroup set, there is a map $E^i:Q^{A-i} \rightarrow Q$ such that
\[
E = \{(x_{A-i},E^i(x_{A-i})) \in Q^{A-i} \times Q\}.
\]
Here, we consider $(x_{A-i},E^i(x_{A-i}))$ as an element in $Q^A$. It is easy to check that $E^i$ is a quasigroup.
\end{remark}

\begin{remark}\label{remark: coordinate partition}
If there is a p-representation of degree $q$ for a matroid $M = (E,r)$ and $S$ is an arbitrary set of size $q$, there is a p-representation $\xi$ of degree $q$ for the matroid whose ground set $\Omega$ is a subset of $S^E$ and the blocks of each partition are \emph{coordinate blocks}, i.e., for each $e \in E$, consider the set $\pi_e(\Omega)$ and for each $a \in \pi_e(\Omega) $ consider the set
\[
\Omega_{e,a} = \{x \in \Omega: \pi_e(x) = a\},
\]
and these are the blocks of the partition $\xi_e$. Moreover, if $f_e:S \rightarrow S$ is a permutation for each $e \in E$, then $f=(f_e)_{e \in E}:S^E \rightarrow S^E$ is a bijection. Coordinate partitions on the ground set $\Omega_1 = f(\Omega)$ is a p-representation of degree $q$ for $M$.
See~\cite[Example~1.2]{Matus1999} and~\cite[Remark~2.7]{KalantariKhazaei2026}.
\end{remark}

\begin{proposition}[Circuits yield quasigroup sets~\cite{Matus1999}]\label{prop: circuits give quasigroups}
Consider a matroid $M = (E,r)$ and a p-representation for $M$ on the ground set $\Omega \subseteq Q^E$ with coordinate partitions. If $C$ is a circuit of $M$ and $2 \leq |C|$, then $\pi_C(\Omega) \subseteq Q^C$ is a quasigroup set.
\end{proposition}

\begin{proposition}[Entropy--partition equivalence~\cite{Matus1999}]\label{prop:entropy-partition}
For every integer $q\ge 2$, a matroid is $q$-entropic if and only if it has a partition representation of degree $q$.
\end{proposition}

\begin{proof}
Suppose $(X_e)_{e\in E}$ is a $q$-entropic representation and let $B$ be a basis. Since every $X_b$ takes at most $q$ values,
\[
H(X_B)=|B|\log q=\sum_{b\in B}H(X_b)
\]
forces the variables $(X_b)_{b\in B}$ to be independent and uniform. Hence $X_B$ is uniform on $q^{r(E)}$ points. Since $H(X_E)=H(X_B)$, the whole tuple is a deterministic function of $X_B$ and is uniform on its support $\Omega$, which has cardinality $q^{r(E)}$.

For $A\subseteq E$, choose a basis $B_A$ of the restriction $M|A$. Since $B_A\subseteq A$ and $H(X_A)=r(A)\log q=H(X_{B_A})$, we have $H(X_A\mid X_{B_A})=0$. Thus $X_A$ is a deterministic function of $X_{B_A}$, and consequently $|X_A(\Omega)|\leq q^{r(A)}$. Its entropy equals the logarithm of this upper bound, so $X_A$ is uniform on exactly $q^{r(A)}$ values. Taking the level partitions of the coordinates gives \Cref{def:partition}.

Conversely, choose $\omega$ uniformly from $\Omega$ and let $X_e$ be the block of $\xi_e$ containing $\omega$. The equal-block condition makes $X_A$ uniform on $q^{r(A)}$ joint values, hence $H(X_A)=r(A)\log q$.
\end{proof}

If $\xi$ is a p-representation for a matroid $M = (E,r)$ with coordinate partitions as in \Cref{remark: coordinate partition} and $f_e : S \rightarrow S$ is an arbitrary bijection (a relabeling) for each $e \in E$, we can use the bijection $f = (f_e)_{e \in E}$ to transform $\xi$ and get a new p-representation with coordinate partitions.

\begin{definition}[Linear partition representations]\label{def: linearization}
Consider a finite field $\F$ and a matroid $M = (E,r)$. Assume $\xi$ is a p-representation of degree $|\F|$ for $M$ on the ground set $\Omega \subseteq \F^E$ with coordinate partitions. If there is a matrix $A \in M_{r(E) \times |E| }(\F)$ such that
\[
\Omega = \{xA : x \in M_{1 \times r(E)}(\F)\},
\]
then $\xi$ is said to be a \emph{linear p-representation} for $M$ over $\F$, and $A$ is called a \emph{generator matrix} for $\xi$.
\end{definition}

Linear p-representations and representations over a finite field (linear representations) are closely related as noted in the following proposition.

\begin{proposition}[Linear partition representations give field representations]\label{prop: linearity}
Consider a finite field $\F$ and a matroid $M = (E,r)$. If $\xi$ is a linear p-representation for $M$ over $\F$ and $A$ is a generator matrix for $\xi$, then $A$ is a representation for $M$ over $\F$.
\end{proposition}

\begin{lemma}[Minor closure of partition representations]\label{lem:minor-closure}
If $M$ has a partition representation of degree $q$, then every minor of $M$ has a partition representation of degree $q$.
\end{lemma}

\begin{proof}
Consider a p-representation of degree $q$ for $M$ on the ground set $\Omega \subseteq Q^E$ with coordinate partitions. If $\omega_0 \in \Omega$, then coordinate partitions on the ground set
\[
\Omega_N = \{\pi_N(x):\pi_J(x) = \pi_J(\omega_0), x\in \Omega\}
\]
gives a degree $q$ representation for the minor $N = M \backslash I / J$.
\end{proof}

\subsection{The foundation of a matroid}

The foundation $F_M$ of a matroid $M$, introduced by
Baker and Lorscheid~\cite{BakerLorscheid2025},
encodes representations of $M$ over fields up to
rescaling. More precisely, for any field $K$,
morphisms from $F_M$ to $K$ correspond to rescaling
classes of $K$-representations of $M$.
The foundation presentation expresses $F_M$ in terms
of the foundations of embedded $U_{2,4}$-minors
and relations among them. We will use this
presentation to construct a representation of the
entire matroid from compatible local data.

The foundation is an algebraic object called a
\emph{pasture}, which generalizes fields by replacing
addition with specified three-term relations.
We recall the definitions needed below.

A \emph{monoid with zero} is a commutative monoid
$P$ with identity $1$ and an absorbing element $0$,
so that $0a=0$ for every $a\in P$.
We denote its group of invertible elements by
$P^\times$.
Let $\operatorname{Sym}_3(P)$ denote the set of
formal sums $a+b+c$, with $a,b,c\in P$, in the
monoid semiring $\mathbb{N}[P]$.

\begin{definition}[Pastures~\cite{BakerLorscheid2025}]
A \emph{pasture} is a monoid with zero $P$ equipped
with a subset $N_P\subseteq\operatorname{Sym}_3(P)$,
called its \emph{nullset}, such that
$P^\times=P\setminus\{0\}$ and, for all $a,b,c,d\in P$,
\begin{enumerate}[label=(P\arabic*)]
\item $a+0+0\in N_P$ if and only if $a=0$;
\item if $a+b+c\in N_P$, then $ad+bd+cd\in N_P$;
\item there is a unique $\epsilon_P\in P^\times$
such that $1+\epsilon_P+0\in N_P$.
\end{enumerate}
We write $-1$ for $\epsilon_P$.
\end{definition}

Every field $K$ is naturally a pasture, with its usual
multiplication and with nullset
\[
N_K=\{a+b+c\in\operatorname{Sym}_3(K):
a+b+c=0\text{ in }K\}.
\]

\begin{definition}[Pasture morphisms~\cite{BakerLorscheid2025}]
A \emph{morphism of pastures} $f:P_1\to P_2$
is a multiplicative map satisfying $f(0)=0$,
$f(1)=1$, and
\[
a+b+c\in N_{P_1}
\quad\Longrightarrow\quad
f(a)+f(b)+f(c)\in N_{P_2}.
\]
\end{definition}

For a pasture $P$ and a set
$S\subseteq\operatorname{Sym}_3(P)$ of relations
$a+b+c$ satisfying $ab\ne 0$, the quotient pasture
$P//S$ is obtained by imposing the relations in $S$
as null relations; see
\cite[Section~2.1.1]{BakerLorscheid2025}.
The following proposition is a special case of
\cite[Proposition~2.6]{BakerLorscheid2025}.

\begin{proposition}[Quotient-pasture factorization~\cite{BakerLorscheid2025}]\label{prop: morphism factors through quotient}
Suppose $P,Q$ are pastures and $S \subseteq \operatorname{Sym}_3(P)$ such that if $a+b+c \in S$, then $ab \neq 0$. If there is a morphism $f:P \rightarrow Q$ such that $f(a) + f(b) +f(c) \in N_Q$ for every $a+b+c \in S$, then there is a morphism $\hat{f}:P//S \rightarrow Q$.
\end{proposition}

If $P_1,\dots,P_n$ are pastures, then it is possible to define the tensor product $\bigotimes_{i=1}^nP_i$, see \cite[Section 2.1.4]{BakerLorscheid2025}. The following proposition is reproduced from \cite[Lemma 2.7]{BakerLorscheid2025}.

\begin{proposition}[Tensor-product universal property~\cite{BakerLorscheid2025}]\label{prop: universal property of tensor product}
Suppose $P_1,P_2,Q$ are pastures. If there are morphisms $f_1:P_1 \rightarrow Q$ and $f_2:P_2 \rightarrow Q$, then there is a morphism $f:P_1 \otimes P_2 \rightarrow Q$ such that $f_1(x) = f(x\otimes 1)$ for all $x \in P_1$ and $f_2(y) = f(1 \otimes y)$ for all $y \in P_2$.
\end{proposition}

By repeated application of the universal property,
morphisms $P_i\to Q$, for $1\leq i\leq m$, induce
a morphism
\[
\bigotimes_{i=1}^{m}P_i\longrightarrow Q.
\]
When a pasture is presented as a quotient of such
a tensor product, the morphism descends to the quotient
if it preserves the defining relations.

To relate foundation morphisms to matroid representations,
we recall Grassmann--Pl\"ucker functions~\cite{BakerLorscheid2025}.

\begin{definition}[Grassmann--Pl\"ucker functions~\cite{BakerLorscheid2025}]\label{def: GP function}
Consider a matroid $M =(E,r)$ and a pasture $P$. Put $r = r(E)$. A function $\Delta:E^r \rightarrow P$ is a \emph{weak Grassmann--Pl\"ucker function} (we use the shorter term \emph{Grassmann--Pl\"ucker function}) if it satisfies the following conditions:
\begin{enumerate}[label=(\roman*)]
\item $\Delta(e_1,\ldots,e_r)\neq 0$ if and only if $\{e_1,\ldots,e_r\}$ is a basis of $M$;
\item $\Delta$ is alternating, i.e., $\Delta(e_{\sigma(1)},\ldots,e_{\sigma(r)})=\operatorname{sign}_P(\sigma)\,\Delta(e_1,\ldots,e_r)$ for every permutation $\sigma$, where $\operatorname{sign}_P(\sigma)=1$ for even $\sigma$ and $\operatorname{sign}_P(\sigma)=\epsilon_P$ for odd $\sigma$;
\item for every $J\in E^{r-2}$ and $e_1,e_2,e_3,e_4\in E$,

\[
\Delta(Je_1e_2)\Delta(Je_3e_4)
+\epsilon_P\Delta(Je_1e_3)\Delta(Je_2e_4)
+\Delta(Je_1e_4)\Delta(Je_2e_3)\in N_P.
\]
\end{enumerate}
\end{definition}

As explained before, a field is an example of a pasture. We have the following proposition about representation of a matroid over a field and a weak Grassmann--Pl\"ucker function for the matroid with values in the field.

\begin{proposition}[Grassmann--Pl\"ucker criterion for representability~\cite{BakerLorscheid2025}]\label{prop: A GP function iff representable over F}
Consider a matroid $M = (E,r)$ and a field $\F$. There is a Grassmann--Pl\"ucker function $\Delta : E^r \rightarrow \F$ if and only if $M$ has a representation over $\F$.
\end{proposition}

\begin{remark}\label{rem: equivalence of represenation}
It will be useful to know how to construct a Grassmann--Pl\"ucker function from a representation of a matroid $M=(E,r)$ over a field $\F$. Note that a representation of a matroid over a field is matrix of size $r \times |E|$ over the field, i.e, a column vector of size $r$ is associated to each member of $E$. The goal is to construct a Grassmann--Pl\"ucker function $\Delta : E^r \rightarrow \F$. Consider a tuple $(e_1,\ldots,e_r) \in E^r$. If $|\{e_1,\ldots,e_r\}| < r$, put $\Delta(e_1,\ldots,e_r) = 0$, otherwise consider the $r\times r$ matrix whose i'th column is the vector associated to $e_i$ and put $\Delta(e_1,\ldots,e_r)$ to be the determinant of this square matrix.
\end{remark}

\begin{definition}[Rescaling equivalence~\cite{BakerLorscheid2025}]\label{def: rescaling equivalent}
Two Grassmann--Pl\"ucker functions $\Delta$, $\Delta^\prime$ with values in the same pasture $P$ are said to be rescaling equivalent if there is a constant $c \in P^\times$ and a map $d:E \rightarrow P^\times$ such that for all $e_1,\ldots,e_r \in E$
\[
\Delta^\prime(e_1,\ldots,e_r) = c . d(e_1)\cdots d(e_r) . \Delta(e_1,\ldots,e_r).
\]
\end{definition}

The foundation of a matroid is a pasture associated to that matroid, see \cite[Definition 4.2]{BakerLorscheid2025}. Consider a matroid $M=(E,r)$. The foundation of $M$ is usually denoted by $F_M$. Put $r = r(M)$. Suppose $J=(j_1,\ldots,j_{r-2}) \in E^{r-2}$ and $\{e_1,e_2,e_3,e_4\} \subseteq E$, with the entries of $J$ and $e_1,e_2,e_3,e_4$ pairwise distinct. Write $\underline{J}=\{j_1,\ldots,j_{r-2}\}$ for the underlying set of $J$. We use the notation $T_B$ of Baker--Lorscheid for the universal Pl\"ucker coordinate associated with a basis $B$. If
\[
\begin{aligned}
\underline{J} &\cup \{e_1,e_3\},\\
\underline{J} &\cup \{e_1,e_4\},\\
\underline{J} &\cup \{e_2,e_4\},\\
\underline{J} &\cup \{e_2,e_3\},
\end{aligned}
\]
are bases, the element
\[
{\crat{e_1}{e_2}{e_3}{e_4}{J}} = \frac{T_{Je_1e_3}.T_{Je_2e_4}}{T_{Je_1e_4}.T_{Je_2e_3}}
\]
is defined. It is called a universal cross ratio of $M$ and it is a member of the foundation $F_M$. Moreover, if $\underline{J} \cup \{e_1,e_2\}$ and $\underline{J} \cup \{e_3,e_4\}$ are also bases, the universal cross ratio is called non-degenerate. Otherwise, it is called degenerate.

\begin{remark}\label{rem:degenerate-cross-ratio}
If an admissible universal cross ratio is degenerate, its value is $1$ in the corresponding foundation. Equivalently, every field-valued Grassmann--Pl\"ucker function evaluates that degenerate cross ratio to $1$; this is the degenerate case of the three-term Pl\"ucker relation.
\end{remark}

The following theorem is from \cite[Theorem 4.6]{BakerLorscheid2025}.
\begin{theorem}[Generators of the foundation~\cite{BakerLorscheid2025}]\label{thm: UCRs as generators}
$F_M^{\times}$ is generated by the collection of all universal ratios.
\end{theorem}

Consider a matroid $M$ and a pasture $P$. Define $\mathcal{X}_M^R(P)$ to be the set of all equivalence classes of Grassmann--Pl\"ucker functions of $M$ with values in $P$. Also, $\operatorname{Hom}(F_M,P)$ is the set of all morphisms from $F_M$ to $P$.
The following theorem is from \cite[Theorem 4.3]{BakerLorscheid2025}.
\begin{theorem}[Foundation morphisms and rescaling classes~\cite{BakerLorscheid2025}]\label{thm: a morphism if and only if a GP function}
Consider a matroid $M$ and a pasture $P$. There is a bijection between  $\mathcal{X}_M^R(P)$ and $\operatorname{Hom}(F_M,P)$.
\end{theorem}

\begin{remark}
The morphism $f:F_M\rightarrow P$ associated with the rescaling class of $\Delta$ by \Cref{thm: a morphism if and only if a GP function} evaluates on a universal cross ratio by
\[
f\left({\crat{e_1}{e_2}{e_3}{e_4}{J}}\right)
=\frac{\Delta(J,e_1,e_3)\Delta(J,e_2,e_4)}
{\Delta(J,e_1,e_4)\Delta(J,e_2,e_3)}.
\]
This is the cross-ratio evaluation formula in \cite[Proposition~3.9]{BakerLorscheid2025}.
\end{remark}
The following theorem is reproduced from~\cite[Corollary~4.21]{BakerLorscheid2025}.
\begin{theorem}[Foundation presentation~\cite{BakerLorscheid2025}]\label{thm: BL foundation representation}
The foundation of a matroid $M = (E,r)$ is naturally isomorphic to a quotient
\[
F_M \simeq \left(\bigotimes_{N}F_N\right)//S,
\]
where $N$ ranges over all embedded minors of $M$ isomorphic to $U_{2,4}$.
The quotient relations are generated by the following four families:
\begin{enumerate}[label=\textup{(\(R\arabic*\))},start=1]
\item\label{rel:R1} If $M\backslash I/J$ is an embedded minor on the ground set $\{e_1,e_2,e_3,e_4,e_5\} = E - (I \cup J)$ and the following universal cross ratios are defined, then
\[
{\crat{e_1}{e_2}{e_3}{e_4}{J}}
.{\crat{e_1}{e_2}{e_4}{e_5}{J}}
.{\crat{e_1}{e_2}{e_5}{e_3}{J}}
=1.
\]
\item\label{rel:R2} If $M\backslash I/J$ is an embedded minor on the ground set $\{e_1,e_2,e_3,e_4,e_5\} = E - (I \cup J)$ and the following universal cross ratios are defined, then
\[
{\crat{e_1}{e_2}{e_3}{e_4}{Je_5}}
.{\crat{e_1}{e_2}{e_4}{e_5}{Je_3}}
.{\crat{e_1}{e_2}{e_5}{e_3}{Je_4}}
=1.
\]
\item\label{rel:R3} If $M\backslash I/J$ is an embedded minor on the ground set $\{e_1,e_2,e_3,e_4,e_5,e_6\} = E - (I \cup J)$, the following universal cross ratios are defined and non-degenerate, and $\cl(Je_5)=\cl(Je_6)$, then
\[
{\crat{e_1}{e_2}{e_3}{e_4}{Je_5}}
=
{\crat{e_1}{e_2}{e_3}{e_4}{Je_6}}.
\]
\item\label{rel:R4} If $M$ has an embedded minor isomorphic to the Fano matroid $F_7$ or its dual $F_7^\ast$, then
\[
-1 = 1.
\]
\end{enumerate}
\end{theorem}

In the quotient-pasture notation, the displayed multiplicative equalities are shorthand for the corresponding three-term null relations (for example, $xyz=1$ is represented by $xyz+(-1)+0\in N$). Thus a morphism from the tensor product descends to $F_M$ exactly when its values satisfy the displayed relations.

\section{Preliminary lemmas}
The foundation presentation reduces the construction of
a global field representation to local representations
of embedded $U_{2,4}$-minors and compatibility
relations among them. We establish the required
local linearizations, study the relevant five-element
configurations, and exclude the Fano matroid and its dual.

\subsection{\texorpdfstring{$U_{2,4}$}{U2,4}}\label{sec: U24}
Embedded $U_{2,4}$-minors provide the local factors
in the foundation presentation. We first use properties
of orthogonal Latin squares of order five to show that
every degree-five partition representation of $U_{2,4}$
can be made linear over $\F_5$ by coordinatewise
relabeling. We then prove that, for a fixed partition
representation, the resulting rescaling class is
independent of the chosen linearization.

\begin{lemma}[Normalization of orthogonal Latin squares~\cite{ColbournDinitz2006}]\label{lem: Cayley latin square}
Consider a Latin square $L:\F_5 \times \F_5 \rightarrow \F_5$ which has an orthogonal mate and $L(0,0) = 0$. There are permutations $f_i:\F_5 \rightarrow \F_5$, $i \in \{1,2,3\}$, such that $f_i(0) = 0$ for each $i$ and $f_3(L(f_1^{-1}(x),f_2^{-1}(y))) = x + y$ for all $x,y \in \F_5$.
\end{lemma}

\begin{proof}
It is well known that if the Latin square $L$ has an orthogonal mate, then there are permutations $g_1,g_2,g_3 : \F_5 \rightarrow \F_5$ such that $g_3(L(g_1^{-1}(x),g_2^{-1}(y))) = x + y$ for all $x,y \in \F_5$, see \cite[Chapter III]{ColbournDinitz2006}. Now consider the permutations $f_1,f_2,f_3:\F_5 \rightarrow \F_5$, where $f_i(x) = g_i(x) - g_i(0)$. Note that $f_i(0) = 0$. By a straightforward computation, it can be shown that for all $x,y \in \F_5$
\[
f_3(L(f_1^{-1}(x),f_2^{-1}(y))) = x + y.
\]
\end{proof}

\begin{lemma}[Orthogonal mates of addition]\label{lem: second needed lem for U24}
Consider a Latin square $R:\F_5 \times \F_5 \rightarrow \F_5$ which is orthogonal to the Latin square $L:\F_5 \times \F_5 \rightarrow \F_5$, $L(x,y) = x+y$. There are a permutation $\pi : \F_5 \rightarrow \F_5$ and $\lambda \in \{2,3,4\}$ such that
\[
R(x,y) = \pi(x+\lambda y).
\]
Moreover, if $R(0,0) = 0$, then $\pi(0) = 0$.
\end{lemma}

\begin{proof}
Fix an arbitrary $t \in \F_5$, and consider the set
\[
R^{-1}(t) = \{(x,y) \in \F_5^2: R(x,y) = t\}.
\]
Since $R$ is a Latin square, there is a permutation $\tau_t: \F_5 \rightarrow \F_5$ such that
\[
R^{-1}(t) = \{(x,\tau_t(x)) \in \F_5^2: x \in \F_5\}.
\]
Since $R$ is orthogonal to $L$, the map $x \mapsto x + \tau_t(x)$ is a permutation. Put $\tilde{\tau}_t: \F_5 \rightarrow \F_5$, $\tilde{\tau}_t(x) = \tau_t(x) - \tau_t(0)$. It is easy to check that $\tilde{\tau}_t$ and the map $x + \tilde{\tau}_t(x)$ are permutations and $\tilde{\tau}_t(0) = 0$. It follows that for each nonzero $x \in \F_5$, $\tilde{\tau}_t(x)$ is nonzero. Put
\[
(a,b,c,d) = (\tilde{\tau}_t(1), \tilde{\tau}_t(2), \tilde{\tau}_t(3), \tilde{\tau}_t(4)).
\]
Note that $a + b  + c + d = 0$, and also
\[
0 = \sum_{x \in \F_5}(x + \tilde{\tau}_t(x))^2 = 2\sum_{x\in\F_5}x\tilde{\tau}_t(x).
\]
Therefore, $a + 2b + 3c + 4d = 0$. It is straightforward to use appropriate linear combinations of $a + b + c + d = 0$ and $a + 2b + 3c + 4d = 0$ to get the relations $a = c + 2d$ and $b = 3c + 2d$. Note that $c,d$ are distinct and nonzero elements of $\F_5$, therefore $c/d \neq 1$. Also, note that $c/d = 3$ gives $a = 0$ and $c/d = 4$ gives $a = d$. We know that $a,d$ are distinct and nonzero elements of $\F_5$, therefore $c/d \neq 3,4$. It follows that $c/d = 2$ and
\[
(a,b,c,d) = (4d,3d,2d,d).
\]
Thus $\tilde{\tau}_t(x) = \lambda_tx$ where $\lambda_t = 4d$. The value $\lambda_t = 4$ is impossible because then $x + \tilde{\tau}_t(x)$ is constant, therefore $\lambda_t \in \{1,2,3\}$.\\
Now, consider two different values $t_1,t_2 \in \F_5$. As we showed,
$R^{-1}(t_1) = \{(x,\lambda_{t_1}x + \tau_{t_1}(0)): x \in \F_5\}$ and
$R^{-1}(t_2) = \{(x,\lambda_{t_2}x + \tau_{t_2}(0)): x \in \F_5\}$. These two sets do not intersect, therefore $\lambda_{t_1} = \lambda_{t_2}$ and $\tau_{t_1}(0) \neq \tau_{t_2}(0)$. Thus there is a $\tilde{\lambda} \in \{1,2,3\}$ such that for every $t \in \F_5$, $R^{-1}(t) = \{(x,\tilde{\lambda}x + \tau_t(0)): x \in \F_5\}$. There is also a permutation $\pi : \F_5 \rightarrow \F_5$ such that $\pi(( -\tilde{\lambda}^{-1})\tau_t(0)) = t$. Putting $\lambda = -\tilde{\lambda}^{-1} \in \{2,3,4\}$ and using the permutation $\pi$ give the required result.\\
Note that if $R(0,0) = 0$, then $\pi(0) = \pi(0+ \lambda \times 0) = R(0,0) = 0$.
\end{proof}

\begin{lemma}[Linearization of $U_{2,4}$]\label{lem: U24 linearity}
Suppose $\xi$ is p-representation of degree $5$ for $U_{2,4}$ on the ground set $\Omega \subseteq \F_5^4$ with coordinate partitions and $(0,0,0,0) \in \Omega$. There are permutations $f_i: \F_5 \rightarrow \F_5$ and $f_i(0) = 0$, for $1 \leq i \leq 4$, such that transforming $\xi$ with $f= (f_1,f_2,f_3,f_4)$ gives the linear p-representation on the ground set
\[
\{(x,y,x + y,x + \lambda y) \in \F_5^4: x,y \in \F_5\},
\]
where $\lambda \in \{2,3,4\} \subseteq \F_5$, with coordinate partitions.
\end{lemma}

\begin{proof}
Since $123$ and $124$ are circuits, there are Latin squares, i.e., binary quasigroups, $L,R:\F_5^2 \rightarrow \F_5$ such that
\[
\Omega = \{(x,y,L(x,y),R(x,y)):x,y \in \F_5\}.
\]
Also, note that $L,R$ are orthogonal mates since $34$ is an independent set in $U_{2,4}$, and $L(0,0) = 0$, $R(0,0) = 0$. By \Cref{lem: Cayley latin square}, there are permutations $f_1,f_2,f_3 : \F_5 \rightarrow \F_5$ such that $f_i(0) = 0$ for each $i$ and putting $f= (f_1,f_2,f_3,\operatorname{id})$, $\Omega_1 = f(\Omega)$, it follows that
\[
\Omega_1 = \{(x,y,x+y,R_1(x,y)):x,y \in \F_5\}.
\]
Here, $R_1$ is a Latin square. Note that $R_1$ is orthogonal to $L_1(x,y) = x+y$ and $R_1(0,0) = 0$.
By \Cref{lem: second needed lem for U24}, there is a permutation $\pi:\F_5 \rightarrow \F_5$ and $\lambda \in \{2,3,4\}$ such that $\pi(0) = 0$ and $R_1(x,y) = \pi(x + \lambda y)$. Now, put $f_4 = \pi^{-1}$ and $g = (f_1,f_2,f_3,f_4)$. The ground set $\Omega_2 = g(\Omega)$ has the required property.
\end{proof}

\begin{lemma}[Uniqueness of the local rescaling class]\label{lem: U24 linearity equivalence}
Let $\Omega$ be a degree-five p-representation of $U_{2,4}$ as in \Cref{lem: U24 linearity}, let $\eta$ be the special linearization supplied by that lemma, and let $\xi$ be any other linearization obtained from the same $\Omega$ by coordinatewise permutations fixing $0$. Then the Grassmann--Pl\"ucker functions associated with generator matrices of $\eta$ and $\xi$ are rescaling equivalent.
\end{lemma}

\begin{proof}
Consider a generator matrix $G$ for $\xi$. Since every two columns of $G$ are independent, $G$ can be chosen to have the following form:
\[G=
\begin{bmatrix}
1 & 0 & a & c\\
0 & 1 & b & d
\end{bmatrix}.
\]
Again, since every two columns of $G$ are independent, $a,b,c,d \in \F_5$ are nonzero. Note that the ground set of $\xi$ is
\[
\Omega_{\xi} = \{(x,y,ax + by, cx + dy): x,y \in \F_5\}.
\]
Also, note that
\[
\Omega_{\eta} = \{(x,y,x+y,x+\lambda y): x,y \in \F_5\}.
\]
Because of $f$, the following equation holds for all $x,y \in \F_5$
\[f_3(x + y) = af_1(x) + bf_2(y). \]
Put $y =0$ in the above equation, and it follows that for every $x$, $f_3(x) = af_1(x)$. Similarly, put $x = 0$, and it follows that for every $y$, $f_3(y) = bf_2(y)$. Therefore, for all $x,y \in \F_5$, it holds that
\[
f_3(x+y) = f_3(x) + f_3(y).
\]
Thus $f_3(x) = xf_3(1)$, and also $f_1(x) = xf_1(1)$ and $f_2(x) = xf_2(1)$.\\
Similarly, the following equation holds for all $x,y \in \F_5$
\[
f_4(x + \lambda y) = cf_1(x) + df_2(y).
\]
Put $y = 0$, it follows that $f_4(x) = cf_1(x)$ for all $x$. Therefore, $f_4(x) = xf_4(1)$. Consequently, multiplying the $i-$th column of a generator matrix for $\eta$ for each $1 \leq i \leq 4$ gives a generator matrix for $\xi$, and the claim follows.
\end{proof}

\subsection{\texorpdfstring{$U_{2,5}$ and $U_{3,5}$}{U2,5 and U3,5}}

Relations (R1) and (R2) in the foundation presentation
involve five-element embedded minors. When all three
universal cross ratios in the respective relation are
non-degenerate, these minors are isomorphic to
$U_{2,5}$ and $U_{3,5}$, respectively.
We establish the coordinatewise linearizations of
degree-five partition representations needed to verify
these relations. The degenerate cases are treated separately.

\begin{lemma}[Linearization of $U_{2,5}$]\label{lem: U25 linearity}
Consider a p-representation $\xi$ of degree five for $U_{2,5}$ such that its ground set $\Omega$ is a subset of $\F_5^5$ with coordinate partitions and $(0,\dots,0) \in \Omega$. There are permutations $f_i:\F_5 \rightarrow \F_5$, for $1 \leq i \leq 5$, such that $f_i(0) = 0$ and transforming $\xi$ with $f = (f_1,\ldots,f_5)$ gives a linear p-representation.
\end{lemma}

\begin{proof}
Since $12$ is a basis of $U_{2,5}$ and $123$, $124$, and $125$ are circuits, there are binary quasigroups (Latin squares) $L_1,L_2,L_3 : \F_5^2 \rightarrow \F_5$ such that
\[
\Omega = \{(x,y,L_1(x,y),L_2(x,y),L_3(x,y)):x,y \in \F_5\}.
\]
Since $34$ is a basis, $L_1$ and $L_2$ are orthogonal mates. By \Cref{lem: Cayley latin square}, there are permutations $g_1,g_2,g_3:\F_5 \rightarrow \F_5$ such that $g_i(0) = 0$, for each $i \in \{1,2,3\}$, and
\[
\Omega_1 = g(\Omega) = \{(x,y,x+y,R_1(x,y),R_2(x,y)):x,y \in \F_5\}.
\]
Here, $g=(g_1,g_2,g_3,g_4,g_5)$, $g_4,g_5 = \operatorname{id}$, and $R_1,R_2$ are Latin squares. Since $34$ and $35$ are bases, by \Cref{lem: second needed lem for U24}, there are permutations $\pi_1,\pi_2 : \F_5 \rightarrow \F_5$ and $\lambda_1,\lambda_2 \in \{2,3,4\}$ such that for all $x,y \in \F_5$
\[
\begin{aligned}
R_1(x,y) &= \pi_1(x + \lambda_1y)\\
R_2(x,y) &= \pi_2(x + \lambda_2y).
\end{aligned}
\]
Moreover, since $R_1(0,0) = R_2(0,0) = 0$, we have $\pi_1(0) = \pi_2(0) = 0$. Put $\pi = (\operatorname{id},\operatorname{id},\operatorname{id},\pi_1^{-1},\pi_2^{-1})$ and $\Omega_2 = \pi(\Omega_1)$. Coordinate partitions on $\Omega_2$ gives a linear p-representation of degree 5 for $U_{2,5}$.
\end{proof}

In order to prove the result about $U_{3,5}$, we need a result from \cite{KokkalaKrotovOstergard2015}. Recall that if $Q$ is a set of size $q$, a subset $C \subseteq Q^n$ is an $(n,M,d)_q$ code if $|C| = M$ and $d = \min_{x\neq y \in C}d_H(x,y)$, where $d_H(x,y)$ is the Hamming distance between $x$ and $y$. Moreover, if $Q$ is a field and $C$ is a subspace of $Q^n$, then $C$ is called linear code. The following theorem is an equivalent formulation of \cite[Theorem III.7]{KokkalaKrotovOstergard2015}.

\begin{theorem}[Linearization of small MDS codes~\cite{KokkalaKrotovOstergard2015}]\label{thm: MDS codes}
Pick $p \in \{5,7\}$. If $k,d \geq 3$ and $C$ is a $(k+d-1,p^k,d)_p$ code, put $n=k+d-1$. Then there are permutations $f_i:\F_p \rightarrow \F_p$, for $1 \leq i \leq n$, such that $f(C)$ is a linear code, where $f=(f_1,\dots,f_n)$.
\end{theorem}

\begin{lemma}[Linearization of $U_{3,5}$]\label{lem: U35 linearity}
Consider a p-representation $\xi$ of degree five for $U_{3,5}$ such that its ground set $\Omega$ is a subset of $\F_5^5$ with coordinate partitions and $(0,\dots,0) \in \Omega$. There are permutations $f_i:\F_5 \rightarrow \F_5$, for $1 \leq i \leq 5$, such that $f_i(0) = 0$ and transforming $\xi$ with $f = (f_1,\ldots,f_5)$ gives a linear p-representation.
\end{lemma}

\begin{proof}
Consider $U_{3,5}$ on the ground set $\{1,2,3,4,5\}$. Since $1234$ and $1235$ are circuits, there are quasigroups $h_1,h_2:\F_5^3 \rightarrow \F_5$ such that
\[
\Omega = \{(a,b,c,h_1(a,b,c),h_2(a,b,c)):a,b,c \in \F_5\}.
\]
Therefore, the Hamming distance between two distinct members of $\Omega$ is at least 3. It is also exactly $3$: fix two of $a,b,c$ and change the third. The changed input coordinate and both quasigroup outputs change, while the other two input coordinates remain fixed. Hence, $\Omega$ is a $(5,5^3,3)_5$ code. By \Cref{thm: MDS codes}, there are permutations $g_i:\F_5 \rightarrow \F_5$, for $1 \leq i \leq 5$, such that $g(\Omega)$ is a linear code, where $g = (g_1, \dots,g_5)$. Write $L=g(\Omega)$ for the resulting linear code and put $f_i(x)=g_i(x)-g_i(0)$ for each $i$. Then $f_i$ is a permutation and $f_i(0)=0$. Moreover, $f(\Omega)=L-(g_1(0),\ldots,g_5(0))$. Since $0\in\Omega$, the vector $(g_1(0),\ldots,g_5(0))$ belongs to $L$, and hence $f(\Omega)=L$ is a linear subspace.
\end{proof}

\subsection{\texorpdfstring{$F_7$ and $F_7^\ast$}{F7 and F7 dual}}\label{sec: fano and its dual}
In the foundation presentation, relation (R4) imposes
$-1=1$ whenever the matroid has an embedded minor
isomorphic to $F_7$ or $F_7^\ast$.
Since this equality does not hold in $\F_5$, we must
exclude these minors.
Degree-five partition representability is closed under
taking minors, so it suffices to show that neither
$F_7$ nor $F_7^\ast$ admits a degree-five
partition representation.

Mat\'u\v{s} proved that $F_7$ admits no degree-five
partition representation
~\cite[Proposition~4.1]{Matus1999}.
We recall his result before establishing the
corresponding exclusion for $F_7^\ast$.

\begin{theorem}[Exclusion of $F_7$~\cite{Matus1999}]\label{thm: fano}
$F_7$ does not have a p-representation of degree 5.
\end{theorem}
Next, we will prove the same result for $F_7^\ast$. In the proof, we will need the following theorem which is the Four Quasigroups Theorem of Shcherbacov~\cite{Shcherbacov2009}.

\begin{theorem}[Four quasigroups theorem~\cite{Shcherbacov2009}]\label{thm: shcherbacov theorem}
Suppose that $Q$ is a non-empty set, and $f_1,f_2,f_3,f_4: Q^2 \rightarrow Q$ are quasigroups. The following equality holds for all $x,y,z \in Q$
\[f_1(f_2(x,y),z) = f_3(x,f_4(y,z))\]
if and only if there is a group structure on $Q$ and there are permutations $\alpha,\beta,\gamma,\mu,\nu: Q \rightarrow Q$ such that
\[
\begin{aligned}
f_1(x,y) &= \mu(x) + \gamma(y)\\
f_2(x,y) &= \mu^{-1}(\alpha(x) + \beta(y))\\
f_3(x,y) &= \alpha(x) + \nu(y)\\
f_4(x,y) &= \nu^{-1}(\beta(x) + \gamma(y)).
\end{aligned}
\]
\end{theorem}
\begin{theorem}[Exclusion of $F_7^\ast$]\label{thm:dual fano}
$F_7^\ast$ does not have a p-representation of degree 5.
\end{theorem}

\begin{proof}
Recall that we consider $F_7^\ast$ to be on the ground set $\{1,\dots,7\}$, and it has the following four-element circuits
\[
1237,\quad 1256,\quad 1346,\quad 1457,\quad 2345,\quad 2467,\quad 3567.
\]
Suppose to the contrary that $F_7^\ast$ has a degree-five p-representation. Consider an arbitrary set $Q$ of size 5. As explained in \Cref{remark: coordinate partition}, we can assume that the ground set of the p-representation $\Omega$ is a subset of $Q^7$. Moreover, the partitions are coordinate partitions. Since $2345$, $1346$, $1237$ are circuits, and $1234$ is a basis, there are quasigroups $h,g,f: Q^3 \rightarrow Q$ such that
\[
\Omega = \{(a,b,c,d,h(b,c,d),g(a,c,d),f(a,b,c)):a,b,c,d \in Q\}.
\]
Fix an arbitrary $b \in Q$ and consider the function $L_b:Q^2 \rightarrow Q$ such that $L_b(c,d) = h(b,c,d)$. Put $\hat{L}_b = \{L_b^{-1}(x):x\in Q\}$. Note that $\hat{L}_b$ is a partition of $Q^2$ with the sets of pre-images of $L_b$, and it has size 5. Fix an arbitrary $a \in Q$ and consider the function $R_a:Q^2 \rightarrow Q$ such that $R_a(c,d) = g(a,c,d)$. As before, define $\hat{R}_a$.

\begin{claim}
For arbitrary $a,b \in Q$, $\hat{L}_b = \hat{R}_a$.
\end{claim}

\begin{proof}
Fix $a,b \in Q$, and consider a block $A \in \hat{L}_b$. Suppose $(c_1,d_1),(c_2,d_2) \in A$, therefore $L_b(c_1,d_1) = L_b(c_2,d_2)$. Note that for $i=1,2$
\[
(a,b,c_i,d_i,L_b(c_i,d_i),R_a(c_i,d_i),f(a,b,c_i)) \in \Omega.
\]
Since $1256$ is a circuit and $\pi_{1256}(\Omega)$ is a quasigroup set, \Cref{prop: circuits give quasigroups} gives $R_a(c_1,d_1)=R_a(c_2,d_2)$. Hence
\[
A\subseteq R_a^{-1}(R_a(c_1,d_1)).
\]
The reverse inclusion follows similarly, so $\hat{L}_b\subseteq\hat{R}_a$. Interchanging the roles of $L$ and $R$ yields the opposite inclusion.
\end{proof}

As proved in the previous claim, there is a partition $\Pi$ of $Q^2$ such that for every $a,b \in Q$, $\Pi = \hat{L}_b = \hat{R}_a$. Next, we want to define a quasigroup $q:Q^2 \rightarrow Q$ such that
\[
\hat{q} = \{q^{-1}(a):a \in Q\} = \Pi.
\]

\begin{claim}
Fix an arbitrary $d_0 \in Q$. There is a binary quasigroup $q:Q^2 \rightarrow Q$ such that $\hat{q} = \Pi$ and $q(c,d_0) = c$ for all $c \in Q$.
\end{claim}

\begin{proof}
Consider the set $S = \{(c,d_0):c\in Q\}$. Note that $|S| = |\Pi| = 5$. Fix an arbitrary $b_0 \in Q$, and consider the function $L_{b_0,d_0}:Q \rightarrow Q$ such that $L_{b_0,d_0}(c) = L_{b_0}(c,d_0) = h(b_0,c,d_0)$. As $h$ is a quasigroup, it follows that $L_{b_0,d_0}$ is a bijection. Moreover, it implies that $|X \cap S| = 1$ for every $X \in \Pi$. Therefore, $S$ is a set of representatives for $\Pi$, and we are able to define a function $q: Q^2 \rightarrow Q$ such that $q(c,d_0) = c$ and $\hat{q} = \Pi$.

It remains to prove that $q$ is a quasigroup. As mentioned, for every $(c,d) \in Q^2$, there is a unique $c_0 \in Q$ such that $L_{b_0}(c,d) = L_{b_0}(c_0,d_0)$, and by definition
\[q(c,d) = L^{-1}_{b_0,d_0}(L_{b_0}(c,d)) = c_0.\]
Now, fix $d \in Q$ and consider $c_1 \neq c_2 \in Q$. Since $L_{b_0}$ is a quasigroup, $L_{b_0}(c_1,d) \neq L_{b_0}(c_2,d)$. Therefore,
\[
q(c_1,d) = L_{b_0,d_0}^{-1}(L_{b_0}(c_1,d)) \neq L_{b_0,d_0}^{-1}(L_{b_0}(c_2,d)) = q(c_2,d).
\]
Similarly, it can be shown that fixing the first variable in $q$ gives a bijection. Therefore, $q$ is a quasigroup.
\end{proof}

Note that $h(b,c,d) = h(b,q(c,d),d_0)$ and $g(a,c,d) = g(a,q(c,d),d_0)$ for all $a,b,c,d \in Q$. Define $A:Q^2 \rightarrow Q$ such that $A(b,c) = h(b,c,d_0)$, and $B:Q^2 \rightarrow Q$ such that $B(a,c) = g(a,c,d_0)$.

\begin{claim}
$A,B$ are quasigroups and
\[
\Omega = \{(a,b,c,d,A(b,q(c,d)),B(a,q(c,d)),f(a,b,c)):a,b,c,d \in Q\}.
\]
\end{claim}

\begin{proof}
Since $h,g$ are quasigroups, $A,B$ are quasigroups as well. Also, note that
\[
\begin{aligned}
A(b,q(c,d)) &= h(b,q(c,d),d_0) = h(b,c,d),\\
B(a,q(c,d)) &= g(a,q(c,d),d_0) = g(a,c,d).
\end{aligned}
\]
Therefore, members of $\Omega$ have the stated form.
\end{proof}
Because $1457$ and $2467$ are circuits, by \Cref{prop: circuits give quasigroups} and \Cref{remark: quasi map eq set}, there are quasigroups $C,D:Q^3 \rightarrow Q$ such that for all $a,b,c,d$
\[
C(a,d,A(b,q(c,d))) = f(a,b,c) = D(b,d,B(a,q(c,d))).
\]
Define $C_1,D_1 : Q^2 \rightarrow Q$ such that $C_1(x,y) = C(x,d_0,y)$ and $D_1(x,y) = D(x,d_0,y)$. Note that for all $a,b,c \in Q$
\[
\begin{aligned}
C_1(a,A(b,c)) &= C_1(a,A(b,q(c,d_0))) = C(a,d_0,A(b,q(c,d_0))) = f(a,b,c)\\
D_1(b,B(a,c)) &= D_1(b,B(a,q(c,d_0))) = D(b,d_0,B(a,q(c,d_0))) = f(a,b,c).
\end{aligned}
\]
Therefore, $C_1(a,A(b,c)) = D_1(b,B(a,c))$ for all $a,b,c \in Q$.

\begin{claim}
$C_1,D_1$ are quasigroups.
\end{claim}

\begin{proof}
By definition, these maps are quasigroups.
\end{proof}

Now, define the quasigroups $A_1,C_2,D_2  : Q^2 \rightarrow Q$ such that for all $x,y,z \in Q$
\[
\begin{aligned}
A_1(x,y) &= z \Longleftrightarrow A(z,y) = x\\
C_2(x,y) &= z \Longleftrightarrow C_1(x,z) = y\\
D_2(x,y) &= z \Longleftrightarrow D_1(z,y) = x.
\end{aligned}
\]
Also, define $C_3:Q^2 \rightarrow Q$ such that $C_3(x,y) = C_2(y,x)$. The fact that $$C_1(a,A(b,c)) = D_1(b,B(a,c))$$ for all $a,b,c \in Q$ is equivalent to $$A_1(C_3(x,y),z) = D_2(x,B(y,z))$$ for all $x,y,z \in Q$. This is the four-quasigroup identity with $(f_1,f_2,f_3,f_4)=(A_1,C_3,D_2,B)$.
By \Cref{thm: shcherbacov theorem}, there is a group structure on $Q$. Since $|Q|=5$, this group is cyclic of order five, and we identify it with $(\F_5,+)$.
There are permutations $\alpha,\beta,\gamma,\mu,\nu: Q \rightarrow Q$ such that
\[
\begin{aligned}
A_1(x,y) &= \mu(x) + \gamma(y)\\
C_3(x,y) &= \mu^{-1}(\alpha(x) + \beta(y))\\
D_2(x,y) &= \alpha(x) + \nu(y)\\
B(x,y) &= \nu^{-1}(\beta(x) + \gamma(y)).
\end{aligned}
\]
Define permutations $f_i:Q \rightarrow Q$, for $1 \leq i \leq 7$, such that
\[
\begin{aligned}
f_1(x) &= -\beta(x), \qquad f_2(x) = x, \qquad f_3(x) =-\gamma(x), \qquad f_4(x) = x,\\ f_5(x) &= \mu(x),
\qquad f_6(x) = -\nu(x), \qquad f_7(x) = \alpha(x).
\end{aligned}
\]
As explained in \Cref{remark: coordinate partition}, coordinate partitions on the set $\Omega_1 = f(\Omega)$ is a p-representation of degree 5 for $F_7^\ast$, where $f=(f_1,\dots,f_7)$. Note that
\[
\Omega_1 = \{(-\beta(a),b,-\gamma(c),d,\mu(A(b,q(c,d))),-\nu(B(a,q(c,d))),\alpha(f(a,b,c))):a,b,c,d \in Q\}.
\]

\begin{claim}
For all $a,b,c,d \in Q$
\[
\begin{aligned}
\mu(A(b,q(c,d))) &= b -\gamma(q(c,d)),\\
\nu(B(a,q(c,d))) &= \beta(a) + \gamma(q(c,d)),\\
\alpha(f(a,b,c)) &= -\beta(a) + b -\gamma(c).
\end{aligned}
\]
\end{claim}

\begin{proof}
Put $X = A(b,q(c,d))$. Therefore,
\[
b = A_1(X,q(c,d)) = \mu(X) + \gamma(q(c,d)),
\]
and the first equation follows. Now, put $Y = B(a,q(c,d))$. It follows that
\[
\nu(Y) = \nu(B(a,q(c,d))) = \nu(\nu^{-1}(\beta(a) + \gamma(q(c,d)))) = \beta(a) + \gamma(q(c,d)),
\]
and the second equation is concluded. Finally, put $Z = f(a,b,c)$. Therefore,
\[
Z = f(a,b,c) = C(a,d_0,A(b,c)) = C_1(a,A(b,c)).
\]
Moreover,
\[
A(b,c) = C_2(a,Z) = C_3(Z,a) = \mu^{-1}(\alpha(Z) + \beta(a)).
\]
It follows that
\[
\alpha(Z) = \mu(A(b,c)) - \beta(a) = b - \gamma(q(c,d_0)) - \beta(a) = b - \gamma(c) - \beta(a),
\]
and the third equation follows.
\end{proof}

\begin{claim}
There is a permutation $\kappa:Q \rightarrow Q$ such that for all $c,d \in Q$
\[
\gamma(q(c,d)) = \gamma(c) + \kappa(d).
\]
\end{claim}

\begin{proof}
Fix $a_1,d_1 \in Q$. By the equations in the previous claim,
\[\alpha(f(a_1,b,c)) = -\beta(a_1) + \mu(A(b,q(c,d_1))) + \gamma(q(c,d_1)) - \gamma(c).\]
Pick an arbitrary $X \in Q$. It is easy to see that for every $c \in Q$, there is a $b_{X,c} \in Q$ such that $X = \mu(A(b_{X,c},q(c,d_1)))$. Therefore,
\[
\alpha(f(a_1,b_{X,c},c)) = -\beta(a_1) + X + \gamma(q(c,d_1)) - \gamma(c).
\]
for all $c \in Q$. Now, the value of the first, fourth, and fifth coordinates are $a_1,d_1$, and $X$, respectively, and the value of the seventh coordinate is $\alpha(f(a_1,b_{X,c},c))$. Since $1457$ is a circuit, it follows that $\gamma(q(c,d_1)) - \gamma(c)$ is independent of $c$. Define $\kappa:Q \rightarrow Q$ such that $\kappa(d) = \gamma(q(c,d)) - \gamma(c)$. It is easy to check that $\kappa$ is a permutation and the claim is proved.
\end{proof}

With the claims that are proved, it follows that
\[
\Omega_1 = \{(-\beta(a),b,-\gamma(c),d,b-\gamma(c)-\kappa(d),-\beta(a) - \gamma(c) -\kappa(d),-\beta(a) + b -\gamma(c)):a,b,c,d \in Q\}.
\]
Now, define $g_4:Q \rightarrow Q$ such that $g_4(d) = -\kappa(d)$, and for each $1 \leq i \leq 7$ and $i \neq 4$, put $g_i = id:Q \rightarrow Q$. Using $g = (g_1,\dots,g_7)$, put $\Omega_2 = g(\Omega_1)$. Note that members of $\Omega_2$ are exactly the same as members of $\Omega_1$, except in the fourth coordinate. Therefore,
\[
\Omega_2 = \{(a,b,c,d,b+c+d,a+c+d,a+b+c):a,b,c,d \in Q\}.
\]
Coordinate partitions on $\Omega_2$ gives a p-representation of degree 5 for $F_7^\ast$. In particular, since $3567$ is a circuit of $F_7^\ast$, $|\pi_{3567}(\Omega_2)| = 5^3$, but it is straightforward to check that $|\pi_{3567}(\Omega_2)| = 5^4$. This is a contradiction.
\end{proof}

\subsection{Structural lemmas}

The verification of relations (R1) and (R2) depends
on which of their three universal cross ratios are
non-degenerate. The following lemmas show that exactly
one non-degenerate cross ratio is impossible, while
three give a uniform minor and two determine a
specific dependent pair or triple. These results
provide the structural information needed to verify
the relations in all cases.

\begin{lemma}[Five-element configurations for (R1)]\label{lem: embedded minor of R1}
Consider a matroid $M  = (E,r)$ and an embedded minor $N = M \backslash I/J$ of $M$. If the ground set of $N$ is $\{e_1,e_2,e_3,e_4,e_5\} = E - (I \cup J)$ and all the cross ratios in (R1) in \Cref{thm: BL foundation representation} are defined, then the following assertions hold:
\begin{enumerate}
\item It is not possible that exactly one of the universal cross ratios in (R1) in \Cref{thm: BL foundation representation} is non-degenerate.
\item If all of the cross ratios in (R1) in \Cref{thm: BL foundation representation} are non-degenerate, then $N$ is isomorphic to $U_{2,5}$.
\item If exactly two of the cross ratios in (R1) in \Cref{thm: BL foundation representation} are non-degenerate, then $N$ has rank 2 and all subsets of the ground set of \(N\) which have size at most 2 are independent except one of its size 2 subsets. More precisely, if the universal cross ratio ${\crat{e_1}{e_2}{e_3}{e_4}{J}}$ (${\crat{e_1}{e_2}{e_4}{e_5}{J}}$) (${\crat{e_1}{e_2}{e_5}{e_3}{J}}$) is degenerate, then $\{e_3,e_4\}$ ($\{e_4,e_5\}$) ($\{e_3,e_5\}$) is dependent.
\end{enumerate}
\end{lemma}

\begin{proof}
Put $A = {\crat{e_1}{e_2}{e_3}{e_4}{J}}$, $B = {\crat{e_1}{e_2}{e_4}{e_5}{J}}$, and $C = {\crat{e_1}{e_2}{e_5}{e_3}{J}}$.
\begin{enumerate}
\item Suppose that the $A$ is non-degenerate and $B$, $C$ are degenerate. Other cases are similar. Since $A$ is non-degenerate, $J \cup \{e_1,e_2\}$ is a basis. This observation with the assumption that $B,C$ are degenerate shows that $J \cup \{e_4,e_5\}$ and $J \cup \{e_3,e_5\}$ are not bases. Also, note that non-degeneracy of $A$ implies that $J \cup \{e_3,e_4\}$ is a basis. Now, since $B$ is defined $J \cup \{e_5\}$ is an independent set and its size is less than $J \cup \{e_3,e_4\}$ which is also an independent set. Therefore, from the properties of independent sets of a matroid, it follows that $J \cup \{e_4,e_5\}$ or $J \cup \{e_3,e_5\}$ is independent, which is a contradiction.
\item Since $J \cup \{e_1,e_4\}$ is a basis, by $A$, $M$ has rank $|J| + 2$, therefore by definition of contraction, it is not possible for $N$ to have a size 3 independent set. Thus $N$ has rank 2. It is straightforward to check that non-degeneracy of $A,B,C$ implies that all two-element subsets of $N$ are independent.
\item Assume that $A$ is degenerate and $B,C$ are non-degenerate. Other cases are similar. Similar to the previous case, $N$ has rank 2. Non-degeneracy of $B$ and $C$ implies 9 two-element subsets of $N$ are independent. The remaining two-element subset is $\{e_3,e_4\}$. In this case, $\{e_3,e_4\}$ is dependent if and only if $A$ is degenerate.
\end{enumerate}
\end{proof}

\begin{lemma}[Five-element configurations for (R2)]\label{lem: embedded minor of R2}
Consider a matroid $M  = (E,r)$ and an embedded minor $N = M \backslash I/J$ of $M$. If the ground set of $N$ is $\{e_1,e_2,e_3,e_4,e_5\} = E - (I \cup J)$ and all the universal cross ratios in (R2) in \Cref{thm: BL foundation representation} are defined, then the following assertions hold:
\begin{enumerate}
\item It is not possible that exactly one of the universal cross ratios in (R2) in \Cref{thm: BL foundation representation} is non-degenerate.
\item If all of the universal cross ratios in (R2) in \Cref{thm: BL foundation representation} are non-degenerate, then $N$ is isomorphic to $U_{3,5}$.
\item If exactly two of the universal cross ratios in (R2) in \Cref{thm: BL foundation representation} are non-degenerate, then $N$ has rank 3 and all subsets of the ground set of \(N\) which have size at most 3 are independent except one of its size 3 subsets. More precisely, if the universal cross ratio ${\crat{e_1}{e_2}{e_3}{e_4}{Je_5}}$ (${\crat{e_1}{e_2}{e_4}{e_5}{Je_3}}$) (${\crat{e_1}{e_2}{e_5}{e_3}{Je_4}}$) is degenerate, then $\{e_1,e_2,e_5\}$ ($\{e_1,e_2,e_3\}$) ($\{e_1,e_2,e_4\}$) is dependent.
\end{enumerate}
\end{lemma}

\begin{proof}
Put $A = {\crat{e_1}{e_2}{e_3}{e_4}{Je_5}} $, $B = {\crat{e_1}{e_2}{e_4}{e_5}{Je_3}}$, and $C = {\crat{e_1}{e_2}{e_5}{e_3}{Je_4}}$.
\begin{enumerate}
\item Suppose that $A$ is non-degenerate and $B,C$ are degenerate. Other cases are similar. Since $A$ is non-degenerate, $J \cup \{e_3,e_4,e_5\}$ is a basis. This observation with the assumption that $B,C$ are degenerate shows that $J \cup \{e_1,e_2,e_3\}$ and $J \cup \{e_1,e_2,e_4\}$ are not bases. Also, note that non-degeneracy of $A$ implies that $J\cup \{e_1,e_2\}$ is an independent set. Now, since $B$ is defined $J \cup \{e_1,e_3,e_4\}$ is a basis and its size is greater than $J\cup \{e_1,e_2\}$, therefore either $J\cup \{e_1,e_2,e_3\}$ or $J\cup \{e_1,e_2,e_4\}$ is independent, which is a contradiction.
\item Since $J \cup \{e_1,e_3,e_5\}$ is a basis, by $A$, $M$ has rank $|J| + 3$, therefore by definition of contraction, it is not possible for $N$ to have a size 4 independent set. Thus $N$ has rank 3. It is straightforward to check that non-degeneracy of $A,B,C$ implies that all three-element subsets of $N$ are independent.
\item Assume that $A$ is degenerate, and $B,C$ are non-degenerate. Other cases are similar. Similar to the previous case, $N$ has rank 3. Non-degeneracy of $B$ and $C$ implies that 9 three-element subsets of $N$ are independent. The remaining three-element subset is $\{e_1,e_2,e_5\}$. In this case, $\{e_1,e_2,e_5\}$ is dependent if and only if $A$ is degenerate.
\end{enumerate}
\end{proof}

\section{Proof of the main theorem}\label{sec:proofmain}

To prove the converse implication in \Cref{thm:main},
we construct a morphism from the foundation of a
$5$-entropic matroid to $\F_5$.
A fixed degree-five partition representation determines
a morphism from the foundation of each embedded
$U_{2,4}$-minor to $\F_5$, independently of the
chosen linearization.
We verify that these local morphisms satisfy the
defining relations of the foundation presentation,
so that they assemble into a global morphism.
This yields the required $\F_5$-representation.

\begin{proof}
First, consider a matroid $M = (E,r)$ that is representable over $\F_5$. We show that $M$ is 5-entropic. Let $(v_e)_{e\in E}$ represent $M$ over $\F_5$ in a finite-dimensional vector space $W$. Choose a linear functional $Y\in W^*$ uniformly and define
\[
X_e=Y(v_e)\in\F_5.
\]
For $A\subseteq E$, the linear map
\[
W^*\longrightarrow\F_5^A,
\qquad
Y\longmapsto(Y(v_e))_{e\in A},
\]
has rank $r(A)$. Its image has $5^{r(A)}$ elements and $X_A$ is uniform on this image. Hence
\[
H(X_A)=r(A)\log5,
\]
so $M$ is 5-entropic.

Conversely, consider a 5-entropic matroid $M = (E,r)$. We show that it is representable over $\F_5$. By \Cref{prop:entropy-partition} and \Cref{remark: coordinate partition}, there is a p-representation $\xi$ of degree 5 for $M$ on the ground set $\Omega \subseteq \F_5^E$ with coordinate partitions. Consider an embedded minor $N = M \backslash I /J$ of type $U_{2,4}$. By \Cref{lem:minor-closure}, coordinate partitions on the ground set
\[
\Omega_N = \{\pi_{E - (I \cup J)}(x) : \pi_J(x) = 0,x \in \Omega \}
\]
gives a p-representation of degree 5 for $N$. From now on, we abuse notation and denote the ground set of $N$, i.e., $E - (I \cup J)$, by $N$ itself. Whenever there are permutations $f_e : \F_5 \rightarrow \F_5$, for each $e \in N$, such that $f_e(0) = 0$ and coordinate partitions on the ground set $\tilde{\Omega}_N = f(\Omega_N)$, where $f = (f_e)_{e \in N}$, gives a linear p-representation for $N$, there is a Grassmann--Pl\"ucker function for $N$ with values in $\F_5$, see \Cref{prop: linearity} and \Cref{rem: equivalence of represenation}. By \Cref{lem: U24 linearity}, there is at least one such family of permutations, which gives a Grassmann--Pl\"ucker function. Moreover, by \Cref{lem: U24 linearity equivalence}, all such Grassmann--Pl\"ucker functions are rescaling equivalent. By \Cref{thm: a morphism if and only if a GP function}, all these rescaling equivalent Grassmann--Pl\"ucker functions give rise to the same morphism from the foundation $F_N$ to $\F_5$.

Up to this point, we have been able to construct a morphism $\Phi_N:F_N \rightarrow \F_5$ for each embedded minor $N$ of $M$ of type $U_{2,4}$. Each $\Phi_N$ is induced from the corresponding minor representation obtained from the common p-representation $\xi$. With these morphisms, using the universal property in \Cref{prop: universal property of tensor product}, it is possible to construct a morphism
\[
\Phi: \left(\bigotimes_{N}F_N\right) \rightarrow \F_5.
\]
For each tensor factor $F_N$, the restriction of $\Phi$ along the canonical factor morphism is precisely $\Phi_N$. We show that $\Phi$ satisfies the relations in \Cref{thm: BL foundation representation}. Coordinate partitions on the ground set
\[
\Omega_K = \{\pi_K(x):\pi_J(x) = (0)_{e \in J}, x\in \Omega\}
\]
gives a p-representation of degree 5 for $K$ and $(0)_{e \in K} \in \Omega_K$.

\begin{claim}
$\Phi$ satisfies relations of type (R1).
\end{claim}

\begin{proof}
Consider an embedded minor $K = M \backslash I / J$ on the ground set $\{e_1,e_2,e_3,e_4,e_5\} = E - (I \cup J)$. We aim to show that $\Phi$ satisfies the relation
\begin{equation}\label{eq : R1 in proof}
{\crat{e_1}{e_2}{e_3}{e_4}{J}}
.{\crat{e_1}{e_2}{e_4}{e_5}{J}}
.{\crat{e_1}{e_2}{e_5}{e_3}{J}}
=1.
\end{equation}
By \Cref{lem: embedded minor of R1}, there are three possibilities with respect to the number of non-degenerate universal cross ratios in the above relation. If all of the universal cross ratios are degenerate, then the relation is $1 = 1$, and $\Phi$ satisfies it trivially. If all of the cross ratios are non-degenerate, then $K$ is isomorphic to $U_{2,5}$, see \Cref{lem: embedded minor of R1}. By \Cref{lem: U25 linearity}, there are permutations $f_i : \F_5 \rightarrow \F_5$, for $1 \leq i \leq 5$, such that $f_i(0) = 0$ and coordinate partitions on the ground set $\Omega_{1,K} = f(\Omega_K)$ gives a linear p-representation for $U_{2,5}$, i.e., there is a matrix $A$ of order $2 \times 5$ with elements in $\F_5$ such that
\[
\Omega_{1,K} = \{xA:x \in \F_5^2\}.
\]
Denote the $i-$th column of $A$ by $A_i$.\\
Consider the universal cross ratio ${\crat{e_1}{e_2}{e_3}{e_4}{J}}$, which is non-degenerate. Therefore, the embedded minor $N = M\backslash (I \cup \{e_5\}) / J$ is of type $U_{2,4}$ on the ground set $\{e_1,e_2,e_3,e_4\}$. Note that transforming the ground set
\[
\Omega_N = \{\pi_N(x):\pi_J(x) = 0,x \in \Omega\}
\]
with $(f_1,f_2,f_3,f_4)$ gives a linear p-representation for $N$ with the generator matrix
\[B = \begin{bmatrix}
A_1 & A_2 & A_3 & A_4
\end{bmatrix},\]
i.e., the first four columns of $A$. By \Cref{prop: linearity}, $B$ gives a representation of $N$ over $\F_5$. We know how to get a Grassmann--Pl\"ucker function for $N$ with values in $\F_5$ given a representation of $N$ over $\F_5$, see \Cref{rem: equivalence of represenation}. By \Cref{thm: a morphism if and only if a GP function}, such a Grassmann--Pl\"ucker function gives a morphism $\Phi_N : F_N \rightarrow \F_5$. By the universal property of tensor product, i.e., \Cref{prop: universal property of tensor product},
\[
\Phi \left( {\crat{e_1}{e_2}{e_3}{e_4}{J}} \right) = \Phi_N \left( {\crat{e_1}{e_2}{e_3}{e_4}{\emptyset}} \right) = \frac{\begin{vmatrix} A_1 & A_3\end{vmatrix}.\begin{vmatrix} A_2 & A_4\end{vmatrix}}{\begin{vmatrix} A_1 & A_4\end{vmatrix}.\begin{vmatrix} A_2 & A_3\end{vmatrix}}.
\]
Here, $|.|$ means determinant, e.g., $\begin{vmatrix} A_1 & A_3\end{vmatrix} = \operatorname{det}\left(\begin{bmatrix} A_1 & A_3\end{bmatrix}\right)$. Similarly,
\[
\begin{aligned}
\Phi \left( {\crat{e_1}{e_2}{e_4}{e_5}{J}} \right) &= \frac{\begin{vmatrix} A_1 & A_4\end{vmatrix}.\begin{vmatrix} A_2 & A_5\end{vmatrix}}{\begin{vmatrix} A_1 & A_5\end{vmatrix}.\begin{vmatrix} A_2 & A_4\end{vmatrix}}\\
\Phi \left( {\crat{e_1}{e_2}{e_5}{e_3}{J}} \right) &= \frac{\begin{vmatrix} A_1 & A_5\end{vmatrix}.\begin{vmatrix} A_2 & A_3\end{vmatrix}}{\begin{vmatrix} A_1 & A_3\end{vmatrix}.\begin{vmatrix} A_2 & A_5\end{vmatrix}}.
\end{aligned}
\]
Therefore,
\[
\Phi \left( {\crat{e_1}{e_2}{e_3}{e_4}{J}} \right) . \Phi \left( {\crat{e_1}{e_2}{e_4}{e_5}{J}} \right).
\Phi \left( {\crat{e_1}{e_2}{e_5}{e_3}{J}} \right) = 1.
\]
So, in the case that all the universal cross ratios in \Cref{eq : R1 in proof} are non-degenerate, the morphism $\Phi$ satisfies the relation. It remains to show that $\Phi$ preserves the relation when two of the universal cross ratios in the relation are non-degenerate.

Assume that ${\crat{e_1}{e_2}{e_4}{e_5}{J}}$ is degenerate and the other two are non-degenerate. Other cases are similar. By \Cref{lem: embedded minor of R1}, all subsets of size at most two of $K$ are independent except $\{e_4,e_5\}$. Recall that coordinate partitions on the ground set
\[
\Omega_K = \{\pi_K(x):\pi_J(x) = (0)_{e \in J}, x \in \Omega\}
\]
gives a degree-five p-representation for $K$. The embedded minor $N = M \backslash (I \cup \{e_5\})/ J$ is isomorphic to $U_{2,4}$. Put
\[
\Omega_N = \{\pi_N(x):\pi_J(x) = (0)_{e \in J}, x \in \Omega\}.
\]
By \Cref{lem: U24 linearity}, there are permutations $f_i: \F_5 \rightarrow \F_5$, for $1 \leq i \leq 4$, and $\lambda \in \{2,3,4\} \subseteq \F_5$ such that $f_i(0) = 0$, for each $i$, and transforming $\Omega_N$ by $f=(f_1,f_2,f_3,f_4)$ gives a linear p-representation
\[
\Omega_{1,N} = f(\Omega_N) = \{(x,y,x+y,x + \lambda y):x,y \in \F_5\}.
\]
Put $g = (f,\operatorname{id})$ and $\Omega_{1,K} = g(\Omega_K)$. Note that
\[
\Omega_{1,N} = \{\pi_N(x):x \in \Omega_{1,K}\}.
\]
Moreover, since $\{e_4,e_5\}$ is dependent and $\{e_4\}$ and $\{e_5\}$ are independent, there is a permutation $\tilde{\pi}$ such that $\tilde{\pi}(0) = 0$ and
\[
\Omega_{1,K} = \{(x,y,x+y,x + \lambda y, \tilde{\pi}(x + \lambda y)):x,y \in \F_5\}.
\]
Transform $\Omega_{1,K}$ by $(\operatorname{id},\operatorname{id},\operatorname{id},\operatorname{id},\tilde{\pi}^{-1})$ to get $\Omega_{2,K}$. Coordinate partitions on the ground set $\Omega_{2,K}$ gives a linear p-representation for $K$ with the generator matrix
\[
G = \begin{bmatrix}
1 & 0 & 1 & 1 & 1\\
0 & 1 & 1 & \lambda & \lambda
\end{bmatrix}.
\]
The matrix
\[
\begin{bmatrix}
1 & 0 & 1 & 1\\
0 & 1 & 1 & \lambda
\end{bmatrix},
\]
i.e., the first four columns of $G$, is a generator matrix for $N$. Therefore,
\[
\Phi \left( {\crat{e_1}{e_2}{e_3}{e_4}{J}} \right) = \Phi_N \left( {\crat{e_1}{e_2}{e_3}{e_4}{\emptyset}} \right) = \frac{\begin{vmatrix}1 & 1\\ 0 & 1\end{vmatrix}.\begin{vmatrix}0 & 1\\ 1 & \lambda\end{vmatrix}}{\begin{vmatrix}1 & 1\\ 0 & \lambda\end{vmatrix}.\begin{vmatrix}0 & 1 \\ 1 & 1\end{vmatrix}} = \frac{1}{\lambda}.
\]
Similarly the embedded minor $M \backslash (I \cup \{e_4\}) / J$ is of type $U_{2,4}$ and coordinate partitions on the ground set
\[
\{(x,y,x+y,x+\lambda y): x,y \in \F_5\}
\]
gives a linear p-representation with generator matrix
\[
\begin{bmatrix}
1 & 0 & 1 & 1\\
0 & 1 & 1 & \lambda
\end{bmatrix}.
\]
Therefore,
\[
\Phi \left( {\crat{e_1}{e_2}{e_5}{e_3}{J}} \right) = \frac{\begin{vmatrix}1 & 1\\ 0 & \lambda\end{vmatrix}.\begin{vmatrix}0 & 1\\ 1 & 1\end{vmatrix}}{\begin{vmatrix}1 & 1\\ 0 & 1\end{vmatrix}.\begin{vmatrix}0 & 1 \\ 1 & \lambda\end{vmatrix}} = \lambda.
\]
Since ${\crat{e_1}{e_2}{e_4}{e_5}{J}}$ is degenerate, \Cref{rem:degenerate-cross-ratio} gives
\[
\Phi\left({\crat{e_1}{e_2}{e_4}{e_5}{J}}\right)=1.
\]
Therefore,
\[
\Phi \left( {\crat{e_1}{e_2}{e_3}{e_4}{J}} \right) . \Phi \left( {\crat{e_1}{e_2}{e_4}{e_5}{J}} \right).
\Phi \left( {\crat{e_1}{e_2}{e_5}{e_3}{J}} \right) = 1.
\]
\end{proof}

\begin{claim}
$\Phi$ satisfies relations of type (R2).
\end{claim}

\begin{proof}
Consider an embedded minor $K = M \backslash I /J$ on the ground set $\{e_1,e_2,e_3,e_4,e_5\} = E - (I \cup J)$. We aim to show that $\Phi$ satisfies the relation
\begin{equation}\label{eq: R2 in proof}
{\crat{e_1}{e_2}{e_3}{e_4}{Je_5}}
.{\crat{e_1}{e_2}{e_4}{e_5}{Je_3}}
.{\crat{e_1}{e_2}{e_5}{e_3}{Je_4}}
=1.
\end{equation}
By \Cref{lem: embedded minor of R2}, there are three possibilities with respect to the number of non-degenerate universal cross ratios in the above relation. If all of the cross ratios are degenerate, then the relation is $1 = 1$, and $\Phi$ satisfies it trivially. If all of the cross ratios are non-degenerate, then $K$ is isomorphic to $U_{3,5}$, see \Cref{lem: embedded minor of R2}. By \Cref{lem: U35 linearity}, there are permutations $f_i:\F_5 \rightarrow \F_5$, for $1 \leq i \leq 5$, such that $f_i(0) = 0$ and coordinate partitions on the ground set $\Omega_{1,K} = f(\Omega_K)$ gives a linear p-representation for $K$, i.e., there is a matrix $A$ of order $3 \times 5$ with elements in $\F_5$ such that
\[
\Omega_{1,K} = \{xA:x \in \F_5^3\}.
\]
Since every three columns of $A$ are independent vectors in $\F_5^3$, $A$ can be chosen so that
\[
A=\begin{bmatrix}
a_1 & a_2 & 1 & 0 & 0\\
b_1 & b_2 & 0 & 1 & 0\\
c_1 & c_2 & 0 & 0 & 1
\end{bmatrix}.
\]
Put $N = M \backslash I / (J \cup \{e_5\})$. Note that coordinate partitions on the ground set
\[
\Omega_{1,N} = \{\pi_N(x):\pi_{J \cup\{e_5\}}(x)=(0)_{e \in J \cup \{e_5\}}, x \in \Omega_{1,K}\}
\]
gives a linear p-representation for $N$ with a generator matrix
\[
\begin{bmatrix}
a_1 & a_2 & 1 & 0\\
b_1 & b_2 & 0 & 1
\end{bmatrix}.
\]
Therefore,
\[
\Phi \left( {\crat{e_1}{e_2}{e_3}{e_4}{Je_5}} \right) = \Phi_N \left( {\crat{e_1}{e_2}{e_3}{e_4}{\emptyset}} \right) = \frac{\begin{vmatrix} a_1 & 1\\ b_1 & 0\end{vmatrix}. \begin{vmatrix}a_2 & 0\\ b_2 & 1\end{vmatrix}}{\begin{vmatrix}a_1 & 0\\ b_1 & 1\end{vmatrix}.\begin{vmatrix}a_2 & 1\\ b_2 & 0\end{vmatrix}}=\frac{(-b_1).a_2}{a_1.(-b_2)}.
\]
Similarly,
\[
\begin{aligned}
\Phi &\left( {\crat{e_1}{e_2}{e_4}{e_5}{Je_3}} \right) = \frac{\begin{vmatrix} b_1 & 1\\ c_1 & 0\end{vmatrix}. \begin{vmatrix}b_2 & 0\\c_2 & 1\end{vmatrix}}{\begin{vmatrix}b_1 & 0\\ c_1 & 1\end{vmatrix}.\begin{vmatrix}b_2 & 1\\ c_2 & 0\end{vmatrix}}=\frac{(-c_1).b_2}{b_1.(-c_2)}\\
\Phi &\left( {\crat{e_1}{e_2}{e_5}{e_3}{Je_4}} \right) = \frac{\begin{vmatrix} a_1 & 0\\ c_1 & 1\end{vmatrix}. \begin{vmatrix}a_2 & 1\\c_2 & 0\end{vmatrix}}{\begin{vmatrix}a_1 & 1\\ c_1 & 0\end{vmatrix}.\begin{vmatrix}a_2 & 0\\ c_2 & 1\end{vmatrix}}=\frac{a_1.(-c_2)}{(-c_1).a_2}.
\end{aligned}
\]
It follows that
\[
\Phi \left( {\crat{e_1}{e_2}{e_3}{e_4}{Je_5}} \right) . \Phi \left( {\crat{e_1}{e_2}{e_4}{e_5}{Je_3}} \right) . \Phi \left( {\crat{e_1}{e_2}{e_5}{e_3}{Je_4}} \right) = 1.
\]
So, $\Phi$ satisfies \Cref{eq: R2 in proof} in the case that all cross ratios are non-degenerate. It remains to show that $\Phi$ preserves the relation in the case that two of the cross ratios are non-degenerate.

Assume that ${\crat{e_1}{e_2}{e_3}{e_4}{Je_5}}$ is degenerate. By \Cref{lem: embedded minor of R2}, all subsets of size at most three of $K$ are independent except $\{e_1,e_2,e_5\}$. Since $\{e_1,e_2,e_3\}$ is independent in $K$, $\{e_1,e_2,e_3,e_4\}$ is a circuit and $\{e_1,e_2,e_5\}$ is a circuit, there are quasigroups $Q:\F_5^3 \rightarrow \F_5$ and $R: \F_5^2 \rightarrow \F_5$ such that
\[
\Omega_K = \{(x,y,z,Q(x,y,z),R(x,y)): x,y,z \in \F_5\}.
\]
Note that $K/e_3$ is isomorphic to $U_{2,4}$. Therefore, for each fixed $z \in  \F_5$, coordinate partitions on the ground set
\[
\{(x,y,Q(x,y,z),R(x,y)):x,y \in \F_5\}
\]
gives a p-representation of degree 5 for $K/e_3$, see \Cref{lem:minor-closure}. In particular, $R$ has an orthogonal mate. It follows that there are permutations $f_1,f_2,f_5 : \F_5 \rightarrow \F_5$ such that $f_1(0) = f_2(0) = f_5(0) = 0$ and transforming $\Omega_K$ by $f = (f_1,f_2,\operatorname{id},\operatorname{id},f_5)$ gives
\[
\Omega_{1,K} = \{(x,y,z,Q_1(x,y,z),x+y):x,y,z \in \F_5\},
\]
where $Q_1$ is a quasigroup. As explained before, for each fixed $z$, the Latin square $(x,y) \mapsto Q(x,y,z)$ is orthogonal to $(x,y) \rightarrow x+y$. Therefore, by \Cref{lem: second needed lem for U24}, there is a permutation $\pi_z : \F_5 \rightarrow \F_5$ and $\lambda_z \in \{2,3,4\} \subseteq \F_5$ such that for all $x,y$
\[
Q_1(x,y,z) = \pi_z(x + \lambda_z y).
\]
Moreover, $\pi_0(0) = 0$. Consider $z_1 \neq z_2 \in \F_5$. If $\lambda_{z_1} \neq \lambda_{z_2}$, then the map
\[
(x,y) \mapsto (x+\lambda_{z_1} y, x + \lambda_{z_2} y)
\]
is a bijection. Thus, there are $x,y \in \F_5$ such that
\[
Q_1(x,y,z_1) = \pi_{z_1}(x+\lambda_{z_1} y) = 0 = \pi_{z_2}(x+\lambda_{z_2} y) = Q_1(x,y,z_2).
\]
Which is a contradiction because $Q_1$ is a quasigroup. Therefore, there is $\lambda \in \{2,3,4\}$ such that for all $z \in \F_5$, $\lambda_z = \lambda$. Now, coordinate partitions on the ground set
\[
\Omega_{N_1} = \{(x,y,\pi_0(x + \lambda y),x+y): x ,y \in \F_5\}
\]
is a p-representation of degree 5 for $N_1 = M \backslash I / (J \cup \{e_3\})$. Transforming $\Omega_{N_1}$ by $(\operatorname{id},\operatorname{id},\pi_0^{-1},\operatorname{id})$ gives
\[
\Omega_{1, N_1} = \{(x,y,x+ \lambda y , x+y):x,y \in \F_5\}.
\]
This p-representation is linear with a generator matrix
\[
\begin{bmatrix}
1 & 0 & 1 & 1\\
0 & 1 & \lambda & 1
\end{bmatrix}.
\]
Therefore,
\[
\Phi \left( {\crat{e_1}{e_2}{e_4}{e_5}{Je_3}} \right) = \Phi_{N_1} \left( {\crat{e_1}{e_2}{e_4}{e_5}{\emptyset}} \right) = \frac{\begin{vmatrix}1 & 1\\ 0 & \lambda\end{vmatrix} . \begin{vmatrix}0 & 1 \\ 1 & 1\end{vmatrix}}{\begin{vmatrix}1 & 1 \\ 0 & 1\end{vmatrix} . \begin{vmatrix}0 & 1 \\ 1 & \lambda\end{vmatrix}} = \lambda.
\]
Put $N_2 = M \backslash I / (J \cup \{e_4\})$. Coordinate partitions on the ground set
\[
\Omega_{N_2} = \{\pi_{N_2}(x): \pi_{e_4}(x) = 0,x \in \Omega_{1,K}\}
\]
gives a p-representation of degree 5 for $N_2$. In fact,
\[
\Omega_{N_2} = \{(\pi_z^{-1}(0) - \lambda y, y ,z, \pi_z^{-1}(0) - \lambda y + y): y,z \in \F_5\}.
\]
Consider $z_1 \neq z_2 \in \F_5$. If $\pi_{z_1}^{-1}(0) = \pi_{z_2}^{-1}(0) = t$, then
\[
Q_1(t,0,z_1) = 0 = Q_1(t,0 ,z_2).
\]
Which is a contradiction because $Q_1$ is a quasigroup. It follows that the map $z \mapsto \pi_z^{-1}(0)$ is a bijection. Also, note that $\pi_0^{-1}(0) = 0$. Denote this bijection by $h$. Transforming $\Omega_{N_2}$ with $(\operatorname{id}, \operatorname{id}, h ,\operatorname{id})$ gives the ground set
\[
\Omega_{1, N_2} = \{(x - \lambda  y,y,x, x - \lambda y + y):x,y \in \F_5\}.
\]
This is a linear p-representation with a generator matrix
\[
\begin{bmatrix}
1 & 0 & 1 & 1\\
-\lambda & 1 & 0 & 1 - \lambda
\end{bmatrix}.
\]
Therefore,
\[
\Phi \left( {\crat{e_1}{e_2}{e_5}{e_3}{Je_4}} \right) = \Phi_{N_2} \left( {\crat{e_1}{e_2}{e_5}{e_3}{\emptyset}} \right) = \frac{\begin{vmatrix}1 & 1\\ -\lambda & 1 - \lambda\end{vmatrix} . \begin{vmatrix}0 & 1 \\ 1 & 0\end{vmatrix}}{\begin{vmatrix}1 & 1 \\ -\lambda & 0\end{vmatrix} . \begin{vmatrix}0 & 1 \\ 1 &1- \lambda\end{vmatrix}} = \frac{1}{\lambda}.
\]
Moreover, since the universal cross ratio ${\crat{e_1}{e_2}{e_3}{e_4}{Je_5}}$ has been assumed to be degenerate, \Cref{rem:degenerate-cross-ratio} gives
\[
\Phi\left({\crat{e_1}{e_2}{e_3}{e_4}{Je_5}}\right)=1.
\]
Therefore,
\[
\Phi \left( {\crat{e_1}{e_2}{e_3}{e_4}{Je_5}} \right)
.\Phi \left( {\crat{e_1}{e_2}{e_4}{e_5}{Je_3}} \right)
.\Phi \left( {\crat{e_1}{e_2}{e_5}{e_3}{Je_4}} \right)
=1.
\]
\end{proof}

\begin{claim}
$\Phi$ satisfies relations of type (R3).
\end{claim}

\begin{proof}
Consider an embedded minor $K = M \backslash I / J$ on the ground set $\{e_1,e_2,\dots,e_6\} = E - (I \cup J)$. We aim to show that $\Phi$ satisfies the relation
\begin{equation}\label{rel: R3 in proof}
{\crat{e_1}{e_2}{e_3}{e_4}{Je_5}} = {\crat{e_1}{e_2}{e_3}{e_4}{Je_6}}.
\end{equation}
As stated in \Cref{thm: BL foundation representation}, both of the universal cross ratios in the above relation are defined and non-degenerate. Moreover, $\operatorname{cl}(Je_5) = \operatorname{cl}(Je_6)$. It follows that $\{e_5,e_6\}$ is a circuit in $K$. Therefore, there is a permutation $f: \F_5 \rightarrow \F_5$ which describes the coordinate of $e_6$ as the value of the $e_5$ coordinate in $\Omega_K$. Moreover, $f(0) =0$. Recall that
\[
\Omega_K = \{\pi_K(x): \pi_J(x) = (0)_{e \in J}, x \in \Omega\},
\]
and coordinate partitions on this ground set gives a p-representation of degree 5 for $K$. Thus,
\[
\{(x,f(x)):x \in \F_5\} = \pi_{\{e_5,e_6\}}(\Omega_K).
\]
Transforming $\Omega_K$ by $(\operatorname{id},\operatorname{id},\operatorname{id},\operatorname{id},\operatorname{id},f^{-1})$ gives the ground set $\Omega_{1,K}$. Put
\[
N_1=K/e_5\backslash e_6,\qquad N_2=K/e_6\backslash e_5.
\]
These are the embedded $U_{2,4}$-minors corresponding to the two cross ratios. Set
\[
\begin{aligned}
\Omega_{N_1}&=\{\pi_{N_1}(x):x_{e_5}=0,\ x\in\Omega_{1,K}\},\\
\Omega_{N_2}&=\{\pi_{N_2}(x):x_{e_6}=0,\ x\in\Omega_{1,K}\}.
\end{aligned}
\]
After the relabeling by $f^{-1}$, the conditions $x_{e_5}=0$ and $x_{e_6}=0$ select the same subset of $\Omega_{1,K}$, and the remaining coordinates are $e_1,\ldots,e_4$. Hence $\Omega_{N_1}=\Omega_{N_2}$. By \Cref{lem: U24 linearity equivalence}, this common labeled partition representation determines the same rescaling class, and hence the same foundation morphism, for the two embedded $U_{2,4}$ occurrences. Therefore,
\[
\Phi\left({\crat{e_1}{e_2}{e_3}{e_4}{Je_5}}\right)
=\Phi_{N_1}\left({\crat{e_1}{e_2}{e_3}{e_4}{\emptyset}}\right)
=\Phi_{N_2}\left({\crat{e_1}{e_2}{e_3}{e_4}{\emptyset}}\right)
=\Phi\left({\crat{e_1}{e_2}{e_3}{e_4}{Je_6}}\right).
\]
\end{proof}

The claims in this proof show that $\Phi$ satisfies our relations (R1), (R2), and (R3). As shown in \Cref{thm: fano} and \Cref{thm:dual fano}, $M$ cannot have an embedded minor isomorphic to $F_7$ or $F_7^\ast$, since it is 5-entropic. Therefore, $\Phi$ factors through the relations in $S$ and there is a morphism
\[
\tilde{\Phi} : \left(\bigotimes_{N}F_N\right)//S \rightarrow \F_5.
\]
By \Cref{thm: BL foundation representation}, \Cref{thm: a morphism if and only if a GP function}, and \Cref{prop: A GP function iff representable over F}, it follows that $M$ is representable over $\F_5$.
\end{proof}

\section{An application to five-symbol ideal secret sharing}\label{sec:ss}
An \emph{access structure} on a finite participant set
$P$ is an upward-closed family $\Gamma\subseteq 2^P$
with $\varnothing\notin\Gamma$.
A perfect secret-sharing scheme for $\Gamma$ consists
of jointly distributed random variables
$\big(S,(X_i)_{i\in P}\big)$, with $H(S)>0$, such that,
for every $A\subseteq P$,
\[
H(S\mid X_A)=
\begin{cases}
0, & A\in\Gamma,\\
H(S), & A\notin\Gamma.
\end{cases}
\]
Here $S$ denotes the secret, $X_i$ is the share of
participant $i$, $X_A=(X_i)_{i\in A}$, and $H$ denotes
Shannon entropy. The scheme is \emph{ideal} if
$H(X_i)=H(S)$ for every $i\in P$.
It is \emph{$q$-ideal} if it is ideal, the secret
is uniformly distributed over a $q$-element alphabet,
and every share takes values in a $q$-element alphabet.

An access structure is called \emph{ideal}
(respectively, \emph{$q$-ideal}) if it admits
an ideal (respectively, $q$-ideal) perfect
secret-sharing scheme.

Let $M=(E,r)$ be a matroid and let $e\in E$ be a nonloop.
The family
\[
\Gamma_e(M)=
\{A\subseteq E\setminus\{e\}:
r(A\cup\{e\})=r(A)\}
\]
is an access structure on $E\setminus\{e\}$, called the
\emph{port} of $M$ at $e$.

Brickell and Davenport showed that every ideal perfect
secret-sharing scheme determines a matroid whose port
at the dealer is the given access structure~\cite{BrickellDavenport1991}.
For a $q$-ideal scheme, the associated matroid is
$q$-entropic. Conversely, a $q$-entropic matroid
determines a $q$-ideal perfect secret-sharing scheme
through any nonloop element chosen as the dealer.
We summarize this correspondence in the following
theorem~\cite{Matus1999,AbbeSpirkl2019}.

\begin{theorem}[Fixed-alphabet secret sharing and matroid ports~\cite{BrickellDavenport1991,Matus1999}]\label{thm:ss-entropic-port}
Let $q\geq 2$ and let $\Gamma$ be an access structure
on a finite participant set $P$. The following are
equivalent:
\begin{enumerate}
\item $\Gamma$ is $q$-ideal.
\item There exists a $q$-entropic matroid $M$ on
$P\cup\{0\}$, with $r_M(\{i\})=1$ for every
$i\in P\cup\{0\}$, such that
\[
\Gamma=\Gamma_0(M).
\]
\end{enumerate}
\end{theorem}

A secret-sharing scheme with secret and shares taking
values in a finite field $K$ is called \emph{$K$-linear}
if there exist a positive integer $d$, a uniformly
distributed random vector $R\in K^d$ (called randomness), and linear
maps $\ell_i:K^d\to K$, for $i\in P\cup\{0\}$,
such that
\[
S=\ell_0(R),\qquad X_i=\ell_i(R)\quad(i\in P).
\]
For such a scheme, a set $A\subseteq P$ is qualified
if and only if
\[
\ell_0\in\operatorname{span}_K
\{\ell_i:i\in A\}.
\]
Thus its access structure is the port at $0$ of
the matroid represented by these linear maps.

\begin{corollary}[Five-symbol ideal secret sharing]\label{cor:ss}
Every $5$-ideal access structure admits an
ideal perfect $\F_5$-linear secret-sharing scheme.
\end{corollary}

Indeed, by \Cref{thm:ss-entropic-port,thm:main},
every $5$-ideal access structure is the port of an
$\F_5$-representable matroid, and therefore admits
an ideal $\F_5$-linear realization
~\cite{BrickellDavenport1991}.

\section{Concluding remarks}\label{sec:conclusion}

We have proved that a finite matroid is $5$-entropic
if and only if it is representable over $\F_5$.
The proof uses the foundation presentation to assemble
compatible local representations into a global one.
The classifications of the five-element degeneracies
are purely matroid-theoretic, and the identification
of the contraction fibres does not depend on the
alphabet size. Their use in constructing a global
field representation, however, requires appropriate
local linearizations and compatibility results.

For another prime alphabet size $q$, the same foundation
presentation reduces the representability question to
constructing $\F_q$-valued morphisms on embedded
$U_{2,4}$-minors and verifying their defining relations.
In particular, whether every $7$-entropic matroid is
representable over $\F_7$ remains open.
Such representability cannot hold for every alphabet
size: the non-Pappus matroid is $9$-entropic but is not
representable over any field
~\cite{SimonisAshikhmin1998}.

\par\medskip\noindent\textbf{Use of generative AI.} Generative-AI tools were used for language editing and exploratory mathematical assistance. The authors independently checked the mathematical statements and take full responsibility for the content.

\end{document}